\documentclass[12pt,leqno,a4paper,twoside]{amsart}
\usepackage{amsmath,amsthm,amssymb,verbatim,enumerate,ifthen,times,cite}
\usepackage[mathscr]{eucal}
\usepackage[utf8]{inputenc}
\usepackage[T1]{fontenc}

\newcommand{\N}{\mathbb{N}}

\newcommand{\A}{\mathscr{A}}

\newcommand{\F}{\mathscr{F}}

\newtheorem{theorem}{Theorem}[section]
\newtheorem*{theorem*}{Theorem}
\def\Thm#1#2{\ifthenelse{\equal{#1}{*}}{\begin{theorem*}#2\end{theorem*}}
  {\begin{theorem}\label{T#1}#2\end{theorem}}}
\newtheorem{Atheorem}{Theorem}

\def\thm#1{Theorem~\ref{T#1}}

\newtheorem{proposition}[theorem]{Proposition}
\newtheorem*{proposition*}{Proposition}
\def\Prp#1#2{\ifthenelse{\equal{#1}{*}}{\begin{proposition*}#2\end{proposition*}}
             {\begin{proposition}\label{P#1}#2\end{proposition}}}
\def\prp#1{Proposition~\ref{P#1}}

\newtheorem{corollary}[theorem]{Corollary}
\newtheorem*{corollary*}{Corollary}
\def\Cor#1#2{\ifthenelse{\equal{#1}{*}}{\begin{corollary*}#2\end{corollary*}}
             {\begin{corollary}\label{C#1}#2\end{corollary}}}

\newtheorem{lemma}[theorem]{Lemma}
\newtheorem*{lemma*}{Lemma}
\def\Lem#1#2{\ifthenelse{\equal{#1}{*}}{\begin{lemma*}#2\end{lemma*}}
             {\begin{lemma}\label{L#1}#2\end{lemma}}}
\def\lem#1{Lemma~\ref{L#1}}

\newtheorem{example}{Example}
\newtheorem*{example*}{Example}
\def\Exa#1#2{\ifthenelse{\equal{#1}{*}}{\begin{example*}\rm #2\end{example*}}
             {\begin{example}\label{Ex#1}\rm #2\end{example}}}

\newtheorem{problem}[theorem]{Problem}

\theoremstyle{definition}
\newtheorem{definition}[theorem]{Definition}
\newtheorem*{definition*}{Definition}
\def\Defi#1#2{\ifthenelse{\equal{#1}{*}}{\begin{definition*}#2\end{definition*}}
             {\begin{definition}\label{D#1}#2\end{definition}}}
\def\defi#1{Definition~\ref{D#1}}

\newtheorem{remark}[theorem]{Remark}
\newtheorem*{remark*}{Remark}
\def\Rem#1#2{\ifthenelse{\equal{#1}{*}}{\begin{remark*}\rm #2\end{remark*}}
             {\begin{remark}\label{R#1}\rm #2\end{remark}}}

\newcommand{\eq}[1]{\eqref{E#1}}
\newcommand{\Eq}[2]{\ifthenelse{\equal{#1}{*}}
  {\begin{equation*}\begin{aligned}[]#2\end{aligned}\end{equation*}}
  {\begin{equation}\begin{aligned}[]\label{E#1}#2\end{aligned}\end{equation}}}

\long\def\comment#1{}

\def\supp{\mathop{\rm supp}}
\def\cl{\mathop{\rm cl}\nolimits}
\def\Cl{\mathop{\rm Cl}\nolimits}
\def\A{\mathscr{A}}

\def\T{\mathscr{T}}

\DeclareMathOperator{\conv}{conv}

\begin{document}

\title[The R{\aa}dström Cancellation Theorem in Cornets]
{An extension of the R{\aa}dström Cancellation Theorem \\ to Cornets}

\author[G. M. Molnár]{Gábor M. Molnár}
\address{Doctoral School of Mathematical and Computational Sciences, University of Debrecen, H-4002 Debrecen, Pf.\ 400, Hungary}
\email{molnar.gabor.marcell@science.unideb.hu}

\author[Zs. P\'ales]{Zsolt P\'ales}
\address{Institute of Mathematics, University of Debrecen, H-4002 Debrecen, Pf.\ 400, Hungary}
\email{pales@science.unideb.hu}

\subjclass[2010]{}

%\dedicatory{}

\thanks{The research of the second author was supported by the K-134191 NKFIH Grant and the 2019-2.1.11-T\'ET-2019-00049, EFOP-3.6.1-16-2016-00022 and EFOP-3.6.2-16-2017-00015 projects. The last two projects are co-financed by the European Union and the European Social Fund.}

\begin{abstract}
The aim of this paper is to introduce the notion of cornets, which form a particular subclass of ordered semigroups also equipped with a multiplication by natural numbers. The most important standard examples for cornets are the families of the nonempty subsets and the nonempty fuzzy subsets of a vector space. In a cornet, the convexity, nonnegativity, Archimedean property, boundedness, closedness of an element can be defined naturally. The basic properties related to these notions are established. The main result extends the Cancellation Principle discovered by R{\aa}dström in 1952.
\end{abstract}
\maketitle

\section{Introduction}

In the theory of convex sets, a basic Cancellation Principle was discovered by R{\aa}dström \cite{Rad52b} in 1952. The Lemma 2 of his paper states that the inclusion 
\Eq{*}{
  A+B\subseteq C+B
}
implies $A\subseteq C$ provided that $A,B,C$ are nonempty subsets of a normed space $X$, $C$ is closed and convex and $B$ is bounded. 

This lemma turned out to be a basic tool in various fields and hundreds of papers have used it by now. For instance, in nonsmooth analysis \cite{BaiFar10,BaiFarRos12a,BaiFarRos12b,CorGajThi10,Dol15,GayGeoMar16,GeoMar18,HuaNin17}, optimization theory \cite{CreKurRoc14,JouSil20,Kur09}, theory of convex sets and functions \cite{BenSim16,CheZho14,DanMedMag11,deBTom11,GrzKucKucUrb14,GrzKucKucUrb15,GrzPalPrzUrb12,GrzPalPrzUrb18,GrzPalUrb10a,GrzPalUrb16,GrzPalUrb19,GrzPrz15,GrzPrz17,GrzPrzUrb13,GrzUrb14,Iva15,Kwi14,VinNag12,VinNag15}, set-valued analysis \cite{AghNouReg15,ChiTre13,Kho19,KurPopRoc15,LeiMerNikSan14,MirMah14,Orl17,Pia09,Pis13}, set-valued differential equations \cite{AzzBou15a,AzzBou15b,BakGab09,GonSha08,Mal15,PloSkr14}, set-valued functional equations \cite{BaiMosPop18,Mai12,Sik15,Sik16,Sik19,Sma09,Szc09,Szc13,Sun17}, iteration theory \cite{AghNou16,AzoGueMatMer10,Gon12,Pis11,SmaSma12,XuNikZha11}, etc.

These applications motivated us to extend the above Cancellation Principle to a more general setting which, possibly, could allow one to apply it to a broader class of problems. It turns out that the natural setting of the Cancellation Principle is a commutative ordered semigroup which is equipped with a multiplication by natural numbers. These structures will be termed \emph{cornets} in our paper. The most important examples for cornets are the families of the nonempty subsets and the nonempty fuzzy subsets of a vector space. In a cornet, one can naturally define the convexity, nonnegativity, Archimedean property, boundedness, closedness of an element. In Sections 2 and 3, we establish the basic properties related to these notions and, finally, in Section 4, we state an abstract form of the Cancellation Principle and also its consequences.

\section{Cornets and convexity properties in cornets}

In the next two definitions, we describe the main structure, the notion of a cornet, that we shall investigate in this paper.

\Defi{1}{
	An ordered triplet $(X,+,\preceq)$ is called an \emph{ordered commutative semigroup} if
	\begin{enumerate}[(i)]
		\item $(X,+)$ is a commutative unital semigroup with a unit element $0$;
		\item $(X,\preceq)$ is a partially ordered set, that is, $\preceq$ is a reflexive, antisymmetric and transitive binary relation on $X$;
		\item For all $x,y,z\in X$ with $x\preceq y$, the inequality $x+z\preceq y+z$ holds.
	\end{enumerate}
If the partially ordered set $(X,\preceq)$ is complete, i.e., every nonempty lower bounded subset of $X$ has a greatest lower bound, then $(X,+,\preceq)$ is called a \emph{complete ordered commutative semigroup}.
}

A unital subsemigroup $(S,+)$ of an ordered commutative semigroup $(X,+,\preceq)$ is obviously an ordered commutative subsemigroup with the ordering restricted to $S$.

In a semigroup $(X,+)$, we naturally have the multiplication by natural numbers which is defined recursively by
\Eq{*}{
  1\cdot x:=x,\qquad (n+1)\cdot x:=n\cdot x+x \qquad(n\in\N).
}
If the semigroup is unital, then we also define $0\cdot x:=0$.
Using induction, one can easily prove that this multiplication obeys the following rules in an ordered commutative semigroup $(X,+,\preceq)$:
	\begin{enumerate}[(i)]
		\item For all $n,m\in\N$ and $x\in X$, $(nm)\cdot x=n\cdot (m\cdot x)$;
		\item For all $n\in\N$ and $x,y\in X$, $n\cdot (x+y)=n\cdot x+n\cdot y$;
		\item For all $n,m\in\N$ and $x\in X$, $(n+m)\cdot x=n\cdot x+m\cdot x$;
		\item For all $n\in\N$ and $x,y\in X$, if $x\preceq y$, then $n \cdot  x\preceq n\cdot y$.
	\end{enumerate}
%In the sequel, we omit the multiplication sign ``$\cdot$'' where it does not cause any confusion.

In the nex definition we present the central concept of our paper.

\Defi{2}{
	An ordered quadruple $(X,+,*,\preceq)$ is called a \emph{cornet} if $(X,+,\preceq)$ is an ordered commutative semigroup and ``$\,*\,$'' is a multiplication of the elements of $X$ by positive integers such that the following conditions hold:
	\begin{enumerate}[(i)]
		%\item For all $n\in\N$, the map $x\mapsto n*x$ is injective;
		\item For all $n,m\in\N$ and $x\in X$, $(nm)*x=n*(m*x)$;
		\item For all $n\in\N$ and $x,y\in X$, $n*(x+y)=n*x+n*y$;
		\item For all $n,m\in\N$ and $x\in X$, $(n+m)*x\preceq n*x+m*x$;
		\item For all $n\in\N$ and $x,y\in X$, the inequality $x\preceq y$ holds if and only if $n * x\preceq n*y$;
		\item $1*x=x$;
		\item $n*0=0$.
	\end{enumerate}
If the partially ordered set $(X,\preceq)$ is complete, then $(X,+,*,\preceq)$ is called a \emph{complete cornet}. 
A unital subsemigroup $(S,+)$ of a cornet $(X,+,*,\preceq)$ which is also closed with respect to the multiplication $*$ is called a \emph{subcornet} of $(X,+,*,\preceq)$ with the ordering restricted to $S$.} 

It is obvious that if $(X,+)$ is a commutative unital semigroup, then $(X,+,\cdot,=)$ is a cornet. The following lemma summarizes the basic properties and connection between the two multiplication operations ``$\cdot$'' and ``$*$''.

\Lem{.}{Let $(X,+,*,\preceq)$ be a cornet. Then the following two assertions hold.
\begin{enumerate}[(i)]
 \item For all $n,k\in\N$, $x_1,\dots,x_k\in X$, 
 \Eq{*}{
  n\cdot(x_1+\dots+x_k)=n\cdot x_1+\dots+n\cdot x_k \quad\mbox{and}\quad
  n*(x_1+\dots+x_k)=n*x_1+\dots+n*x_k.
 } 
 In particular, for all $n,m\in\N$ and $x\in X$, 
 \Eq{nmx}{
   n*(m\cdot x)=m\cdot(n*x).
 }
 \item For all $n,k_1,\dots,k_n\in\N$ and $x\in X$, 
 \Eq{*}{
   (k_1+\cdots+k_n)\cdot x=k_1\cdot x+\dots+k_n\cdot x \quad\mbox{and}\quad
   (k_1+\cdots+k_n)*x\preceq k_1*x+\dots+k_n*x.
 }
 In particular, for all $n,m\in\N$ and $x\in X$,
 \Eq{nx}{
   (mn)*x\preceq n\cdot (m*x).
 }
\end{enumerate}
}

\begin{proof} We prove (i) by induction on $k$. If $k=1$, then the equalities hold trivially. The $k=2$ case follows from property (ii) of the two operations ``$\cdot$'' and ``$*$''. Assume that (i) holds for some $k\in\N$ and let $n\in\N$ and $x_1,\dots,x_{k+1}\in X$ be arbitrary. Then, by property (ii) of the operation ``$*$'' and the inductive hypothesis, we get
\Eq{*}{
  n*(x_1+\dots+x_k+x_{k+1})
  =n*(x_1+\dots+x_k)+n*x_{k+1}
  =n*x_1+\dots+n*x_k+n*x_{k+1}.
}
For the operation ``$\cdot$'', the proof is completely similar. 

By taking $k:=m$ and $x_1:=\dots =x_k:=x$, the second equality in (i) yields the equality \eq{nmx}.

The relations in (ii) will be proved by induction on $n$. For $n=1$ both of them hold with equality. For $n=2$, they are consequences of property (iii) of the two operations ``$\cdot$'' and ``$*$''. Assume that (ii) holds for some $n\in\N$ and let $k_1,\dots,k_{n+1}\in\N$ and $x\in C$ be arbitrary. Then, by property (iii) of the operation ``$*$'' and the inductive hypothesis, we get
\Eq{*}{
  (k_1+\cdots+k_n+k_{n+1})*x
  \preceq(k_1+\cdots+k_n)*x+k_{n+1}*x
  \preceq k_1*x+\dots+k_n*x+k_{n+1}*x.
}
For the operation ``$\cdot$'', the proof is completely similar.

By taking $k_1:=\dots =k_n:=m$, the second inequality in (ii) yields property \eq{nx}.
\end{proof}

For a given element $x\in X$, the set of those numbers $n$ for which \eq{nx} holds with equality if $m=1$ play a crucial role among the properties of $x$.

\Defi{2,5}{Let $(X,+,*,\preceq)$ be a cornet and $n\in\N$. An  element $x\in X$ will be called \emph{$n$-convex} if it fulfills the equality $n*x=n\cdot x$. For fixed elements $x\in X$ and $n\in\N$, we introduce the notations
\Eq{*}{
  C_x:=\{n\in\N\mid x\mbox{ is $n$-convex}\} \qquad\mbox{and}\qquad
  C^n:=\{x\in X\mid x\mbox{ is $n$-convex}\},
}
respectively. If $C_x=\N$, i.e., if $x$ is $n$-convex for all $n\in\N$, then we say that $x$ is \emph{convex}.}

\Lem{C}{Let $(X,+,*,\preceq)$ be a cornet. Then the following assertions hold:
\begin{enumerate}[(i)]
 \item For all $x\in X$, the set $C_x$ is a unital multiplicative subsemigroup of $\N$.
 \item For all $n\in\N$, the quadruple $(C^n,+,*,\preceq)$ is a subcornet of $(X,+,*,\preceq)$.
\end{enumerate}
}

\begin{proof} Let $x\in X$ be fixed. It is clear that $1\in C_x$. Let $n,m\in C_x$ be arbitrary. Then, by property (i) of the two multiplication operations, by the $n$- and $m$-convexity of $x$ and by \eq{nx}, we have that
\Eq{*}{
   (mn)*x=m*(n*x)=m*(n\cdot x)=n\cdot(m* x)=n\cdot(m\cdot x)=(nm)\cdot x.
}
This shows that $x$ is also $(mn)$-convex, i.e., $mn\in C_x$.

For the second assertion, let $n\in\N$ be fixed and $x,y\in C^n$. Using property (ii) of the two multiplication operations and the $n$-convexity of $x$ and $y$, we have
\Eq{*}{
  n*(x+y)=n*x+n*y=n\cdot x+n\cdot y=n\cdot(x+y),
}
therefore, $x+y$ is also $n$-convex.

If $x\in C^n$ and $m\in\N$, then
\Eq{*}{
  n*(m*x)=(nm)*x=(mn)*x=m*(n*x)=m*(n\cdot x)=n\cdot(m*x),
}
which proves that $m*x$ is $n$-convex.  
\end{proof}

In what follows, we define the $n$-convex hull of elements in a cornet $(X,+,*,\preceq)$.

\Defi{ch}{Let $n\in\N$, $(X,+,*,\preceq)$ and $x\in X$. The \emph{$n$-convex hull} of $x$, denoted as $\conv_n(x)$, is the smallest element $y\in C^n$ such that $x\preceq y$, that is, whenever $x\preceq z\in C^n$, then $y\preceq z$.}

In general, the $n$-convex hull of an element may not exist. In order to formulate conditions which are sufficient for the existence, we say that the $*$-multiplication in a complete cornet $(X,+,*,\preceq)$ is $n$-continuous (with respect to the ordering "$\preceq$") if, for all nonempty lower bounded subsets $H\subseteq X$, we have
\Eq{*}{
   \inf(n*H)=n*\inf(H).
}

\Prp{0}{Let $n\in\N$ and let $(X,+,*,\preceq)$ be a complete cornet in which the $*$-multiplication is $n$-continuous. Then $(C^n,+,*,\preceq)$ is a complete  subcornet of $(X,+,*,\preceq)$. Furthermore, for every element $x\in X$, $x$ admits an $n$-convex hull if and only if it has an $n$-convex majorant.}

\begin{proof} Let $H\subseteq C^n$ be a lower bounded subset and denote $x:=\inf(H)$. Then, for all $h\in H$, 
\Eq{*}{
   n\cdot x\preceq n\cdot h=n*h,
}
hence
\Eq{*}{
   n\cdot x\preceq\inf(n*H)=n*\inf(H)=n*x.
}
The reversed inequality is a consequence of \eq{nx} with $m=1$, hence $n\cdot x=n*x$ holds, which shows that $x$ is also $n$-convex. This proves that $(C^n,\preceq)$ is a complete partially ordered set. 

To prove the last assertion, let $x\in X$ be arbitrary. If $x$ has an $n$-convex hull, then it also has an $n$-convex majorant. Conversely, if $x$ admits an $n$-convex majorant, then the set 
\Eq{*}{
 H:=\{z \in C^n\mid x\preceq z\}
}
is nonempty and lower bounded. According to the first part, the infimum $u$ of $H$ belongs to $C^n$, that is, $u$ is $n$-convex. 
It is clear that $u$ is the $n$-convex hull of $x$.
\end{proof}

In a cornet $(X,+,*,\preceq)$, let $K^n$ denote the collection of those elements which have an $n$-convex hull and let $M^n$ denote the set of those elements that have an $n$-convex majorant. Obviously, we have $C^n\subseteq K^n\subseteq M^n$. Using this terminology, the previous proposition asserts that if $(X,+,*,\preceq)$ is a complete cornet in which the $*$-multiplication is $n$-continuous, then $K^n=M^n$.

\Prp{co}{Let $(X,+,*,\preceq)$ be a cornet and let $n\in N$. Then we have the following assertions.
	\begin{enumerate}[(i)]
		\item If $x\in K^n$, then $\conv_n(x)\in C^n$ and $x\preceq\conv_n(x)$. Furthermore, $\conv_n:K^n\to C^n$ is a monotone mapping whose set of fixed points is equal to $C^n$.
		\item $(M^n,+,*,\preceq)$ is a subcornet of $(X,+,*,\preceq)$. Furthermore,
		\Eq{cv+}{
		\conv_n(x+y)&\preceq\conv_n(x)+\conv_n(y) 
		   \qquad&&\mbox{if } x,y,x+y\in K^n, \\
		\conv_n(m*x)&\preceq m*\conv_n(x)\qquad&&\mbox{if } x,m*x\in K^n.
		}
	\end{enumerate}
}

\begin{proof}\hfill\null

	(i) For an arbitrary $x\in K^n$, the inclusion $\conv_n(x)\in C^n$ and the inequality $x\preceq\conv_n(x)$ are consequences of the definition of the $n$-convex hull. If $x\in C^n$, then the smallest $n$-convex element which is nonsmaller than $x$ is equal to $x$, that is, $x=\conv_n(x)$. Conversely, if $x=\conv_n(x)$, then $\conv_n(x)\in C^n$ implies that $x$ must be in $C^n$.
	To see that $\conv_n$ is monotone, let $x,y\in K^n$ with $x\preceq y$. Then $x\preceq\conv_n(y)$, which yields that $\conv_n(x)\preceq\conv_n(y)$.
	
	(ii) Let $x,y\in M^n$. Then there exist $u,v\in C^n$ such that $x\preceq u$ and $y\preceq v$. Thus yields that $x+y\preceq u+v\in C^n$, which proves that $x+y\in M^n$. If additionally $x,y,x+y\in K^n$, then the inequalities $x\preceq\conv_n(x)$ and $y\preceq\conv_n(y)$ imply that $x+y\preceq\conv_n(x)+\conv_n(y)\in C^n$. This proves the first inequality in \eq{cv+}.
	
	For the second inequality in \eq{cv+}, let $x\in M^n$ and $m\in\N$. Then there exist $u\in C^n$ such that $x\preceq u$. This implies $m*x\preceq m*u\in C^n$, which shows that $m*x\in M^n$. If additionally $x,m*x\in K^n$, then the inequality $x\preceq\conv_n(x)$ yields that $m*x\preceq m*\conv_n(x)\in C^n$. This shows the second assertion of \eq{cv+}.
\end{proof}

To illustrate the rich applicability of the above concepts, we provide the most basic examples for cornets in the subsequent three propositions.
For these definitions, we introduce the notion of wedge in abelian group setting. 

\Defi{W}{If $(G,+)$ is an abelian semigroup and $n\in\N$, then for a subset $S\subseteq G$, define
\Eq{*}{
  n^{-1}(S):=\{x\in G\mid n\cdot x\in S\}.
}
A subsemigroup $S$ of the group $(G,+)$ is said to be \textit{$n$-divisible} if, for all $x\in S$, the set $n^{-1}(\{x\})\cap S$ is nonempty. If this set is a singleton, then $S$ is called \textit{uniquely $n$-divisible} and its unique element will be denoted by $x/n$. \\\indent In a unital abelian semigroup $G$, a subset $W\subseteq G$ is called a \emph{wedge} if the following properties are satisfied:
\begin{enumerate}[(i)]
 \item $W$ is a unital subsemigroup of $G$.
 \item If $u,v\in W$ such that $u+v=0$, then $u=v=0$.
 \item For all $n\in\N$, the inverse image $n^{-1}(W)$ is contained in $W$.
\end{enumerate}
In terms of a wedge $W\subseteq G$, we can define a partial order $\preceq_W$ in the following way: For $x,y\in G$, we say that $x\preceq_W y$ if $y\in x+W$. It immediately follows that $\preceq_W$ is a reflexive, and transitive relation on $G$. If, in addition, $G$ is cancellative (which is always the case if $G$ is group), then $\preceq_W$ is antisymmetric and hence it is a partial order on $G$.}

\Prp{1+}{Let $(G,+)$ be a abelian group and let $W\subseteq G$ be a wedge. Then, for a subsemigroup $S$ of $G$ containing $W$, the quadruple $(S,+,\cdot,\preceq_W)$ is a cornet in which every element is $n$-convex for all $n\in\N$. In particular, by taking $W:=\{0\}$, it follows that $(G,+,\cdot,=)$ is a cornet.}

\begin{proof}
The properties (i), (ii), (iii) of \defi{2} can easily be verified by induction, moreover, (iii) holds with equality. Thus, it suffices to show that property (iv) is also valid. 

Let $n\in\N$ and $x,y\in S$ be arbitrary. Assume first that $x\preceq_W y$ holds. Then $y\in x+W$. The set $W$ is a subsemigroup, therefore, $y\in x+W$ implies that $n\cdot y\in n\cdot x+n\cdot W\subseteq n\cdot x+W$, which yields that $n\cdot x\preceq_W n\cdot y$. On the other hand, if $n\cdot x\preceq_W n\cdot y$ holds, then $n\cdot(y-x)\in W$, consequently $y-x\in n^{-1}(W)$. By condition (iii) of \defi{W}, it follows that $y-x\in W$ must be valid and hence $x\preceq_W y$.

The operation $\cdot$ being the cornet-multiplication implies that every element of $S$ is $n$-convex for all $n\in\N$.
\end{proof}

\Prp{2}{Let $(G,+)$ be an abelian group, $W$ be a wedge and let $S$ be a subsemigroup of $G$ containing $W$. Let $P_W(S)$ denote the collection of all nonempty $W$-invariant subsets $A$ of $S$, which means that $A+W\subseteq A$ holds. Define the operations $+$ and $*$ by:
\Eq{+*}{
  A+B&:=\{a+b\mid a\in A,\,b\in B\} &&\qquad(A,B\in P_W(S)),\\
  n*A&:=\{n\cdot a+w\mid a\in A,\,w\in W\} &&\qquad
  (A\in P_W(S),\,n\in\N).
}
Then $(P_W(S),+,*,\subseteq)$ is a complete cornet with the unit element $W$. Furthermore, the mapping
\Eq{vph}{
  \varphi(x):=x+W \qquad(x\in S)
}
is an injective order reversing homomorphic mapping of $(S,+,\cdot,\preceq_W) $ into $(P_W(S),+,*,\subseteq)$. In addition, if $n\in\N$ and $W$ is $n$-divisible, then $A\in P_W(S)$ is $n$-convex if and only if, for all $x_1,\dots,x_n\in A$, we have 
\Eq{n-1}{
n^{-1}(\{x_1+\dots+x_n\})\cap A\neq\emptyset.}
}

\begin{proof} 
If $A,B\in P_W(S)$, then $A+B\subseteq S+S\subseteq S$ and $(A+B)+W=A+(B+W)\subseteq A+B$, which show that $A+B\in P_W(S)$. Therefore, $P_W(S)$ is an abelian semigroup with the addition defined in \eq{+*}. Clearly, for $A\in P_W(S)$, the property $0\in W$ implies $A\subseteq A+W\subseteq A$, which proves that $W$ is the unit element of the semigroup $(P_W(S),+)$.

The inclusion of sets is trivially a partial order on $P_W(S)$ and the implication $A\subseteq B\Rightarrow A+C\subseteq B+C$ is also obvious for $A,B,C\in P_W(S)$. Therefore, $(P_W(S),+,\subseteq)$ is an ordered abelian semigroup.

First observe that the definition of the multiplication operation $*$ is correct, i.e., $n*A\in P_W(S)$ for all $n\in\N$ and $A\in P_W(S)$.

To see that property (i) holds, let $n,m\in\N$ and $A\in P_W(S)$. First, let $u\in (nm)*A$. Then there exist elements $a\in A$ and $w\in W$ such that $u=(nm)\cdot a+w$. We have that $m\cdot a\in m*A$. Therefore,
\Eq{*}{
  u=n\cdot(m\cdot a)+w\subseteq n*(m*A).
}
This proves that $(nm)*A\subseteq n*(m*A)$. To verify the reversed inclusion, let $u\in n*(m*A)$. Then there exist $b\in m*A$ and $w\in W$ such that $u=n\cdot b+w$. Similarly, there exist $a\in A$ and $z\in W$ such that $b=m\cdot a+z$. Combining these equalities, we get that 
\Eq{*}{
u=n\cdot(m\cdot a+z)+w=(nm)\cdot a+(n\cdot z+w)\in (nm)*A,
}
which completes the proof of the reversed inclusion $n*(m*A)\subseteq (nm)*A$ and property (i) of \defi{2}.

The verification of property (ii) of \defi{2} is similar, and therefore it is left to the reader.

To show that (iii) of \defi{2} holds, let $n,m\in\N$ and $A\in P_W(S)$ and $u\in(n+m)*A$. Then there exist $a\in A$ and $w\in W$ such that $u=(n+m)\cdot a+w=n\cdot a+(m\cdot a+w)\in n*A+m*A$. Therefore, $(n+m)*A\subseteq n*A+m*A$. 

For the proof of property (iv) of \defi{2}, let $n\in\N$, $A,B\in P_W(S)$. If $A\subseteq B$ holds, then the inclusion $n*A\subseteq n*B$ is obvious.
Conversely, assume that $n*A\subseteq n*B$ holds. Then, for an arbitrary $a\in A$, we get that $n\cdot a\in n*A\subseteq n*B$, therefore, there exist $b\in B$ and $w\in W$ such that $n\cdot a=n\cdot b+w$. This yields that $n\cdot(a-b)\in W$, i.e., $a-b\in n^{-1}(W)$. Now the condition (iii) of \defi{W} gives that $a-b\in W$, which proves that $a=b+w\in B$. 

The properties (v) and (vi) of \defi{2} can easily be seen.

We verify now the completeness of $(P_W(S),+,*,\subseteq)$. Let $\A:=\{A_\gamma\mid\gamma\in\Gamma\}$ be a nonempty and lower bounded family of elements of $P_W(S)$. Define $A\subseteq S$ by $A:=\bigcap_{\gamma\in\Gamma}A_\gamma$. Then $A$ is nonempty by the existence of a lower bound for the family $\A$.
We show that $A$ also belongs to $P_W(S)$. Indeed, if $u\in A+W$, then there exists $a\in A$ and $w\in W$ such that $u=a+w$. We have that $a\in A_\gamma$ for all $\gamma\in\Gamma$, therefore, $u=a+w\in A_\gamma+W\subseteq A_\gamma$ for all $\gamma\in\Gamma$. This proves $u+w\in A$ showing that $A+W\subseteq A$ holds. Thus, $A\in P_W(S)$ is valid. Clearly, $A$ is the greatest lower bound for $\A$ in $P_W(S)$, and hence $(P_W(S),\subseteq)$ is a complete partially ordered set.

Now consider the mapping $\varphi$ defined by \eq{vph}. For $x\in S$, we have that $\varphi(x)=x+W\subseteq S$ and $\varphi(x)+W=x+W+W\subseteq x+W$, which show that $\varphi(x)$ is in $P_W(S)$. If, for some $x,y\in S$, the equality $\varphi(x)=\varphi(y)$ holds, then $x+W=y+W$, which yields that $x\in y+W$ and $y\in x+W$. Therefore, $x-y\in W\cap(-W)=\{0\}$, showing that $x=y$, which proves the injectivity of $\varphi$.

The structure preserving properties are easily seen from the following identities.
\Eq{*}{
  \varphi(x+y)&=x+y+W=(x+W)+(y+W)=\varphi(x)+\varphi(y), \\
  \varphi(n\cdot x)&=n\cdot x+W=n\cdot(x+W)+W
  =n\cdot\varphi(x)+W=n*\varphi(x).
}
To see that $\varphi$ is also order reversing, observe that
the following inequalities of inclusions are pairwise equivalent:
\Eq{*}{
  x\preceq_W y \quad\Leftrightarrow\quad
  y\in x+W \quad\Leftrightarrow\quad
  y+W\subseteq x+W \quad\Leftrightarrow\quad
  \varphi(y)\subseteq\varphi(x). 
}
Therefore, $\varphi$ is an order reversing and homomorphic embedding of $(S,+,\cdot,\preceq_W) $ into $(P_W(S),+,*,\subseteq)$.

Finally, let $n\in\N$, $A\in P_W(S)$ and assume that $W$ is $n$-divisible. First suppose that $A$ is $n$-convex. Let $x_1,\dots,x_n\in A$ be arbitrary. Then 
\Eq{*}{
x_1+\dots+x_n\in n\cdot A\subseteq n*A=\{n\cdot a+w\mid a\in A,\, w\in W\}.
}
Therefore, there exist $a\in A$ and $w\in W$ such that $x_1+\dots+x_n=n\cdot a+w$.
By the $n$-divisibility of $W$, $n^{-1}(\{w\})\cap W$ is nonempty, therefore, $w=n\cdot v$ for some $v\in W$. On the other hand, $a+v\in A+W\subseteq A$, thus $a+v\in n^{-1}(\{x_1+\dots+x_n\})\cap A$.

To prove the converse, assume that \eq{n-1} is valid for all $x_1,\dots,x_n\in A$. To prove that $A$ is $n$-convex, it is sufficient to show that $n\cdot A\subseteq n*A$.
Let $x\in n\cdot A$. This means that there exist $x_1,\dots,x_n\in A$ such that $x=x_1+\dots+x_n$. By \eq{n-1}, for some $a\in A$, we have that $x_1+\dots+x_n=n\cdot a\in n*A$. Thus $x\in n*A$, which was to be proved.
\end{proof}

\Prp{3}{Let $(G,+)$ be an abelian group, $W$ be a wedge and let $S$ be a  uniquely divisible subsemigroup of $G$ which contains $W$. Let, for $p\in\,]0,1]$,
\Eq{Fp}{
		F^p_W(S):=\{f:S\to[0,1]\mid \sup f\geq p\, \mbox{ and $f$ is $W$-nondecreasing}\}
}
and define the addition and the scalar multiplication in $F_W^p(S)$ by
\Eq{oper}{
  (f\oplus g)(x)
  :=\!\sup_{\substack{u,v\in S\\ u+v=x}}\! \min\big(f(u),g(v)\big),\quad
  (n\odot f)(x):=f\left(\frac{x}{n}\right) 
  \quad(f,g\in F_W^p(S),\,x\in S,\,n\in\N).
}
Finally, let $\leq$ denote the pointwise ordering in $F_W^p(S)$.
Then $(F_W^p(S),\oplus,\odot,\leq)$ is a complete cornet whose unit element is the characteristic function of the wedge $W$. Furthermore, the mapping
\Eq{Phi}{
  \Phi(A):=\chi_A \qquad(A\in P_W(S))
}
is an injective cornet-preserving mapping of $(P_W(S),+,*,\subseteq)$ into $(F_W^1(S),\oplus,\odot,\leq)$. In addition, a function $f\in F_W^p(S)$ is $n$-convex if and only if it is $n$-quasiconcave, i.e., for all $x_1,\dots,x_n\in S$,
\Eq{qc}{
  \min(f(x_1),\dots,f(x_n))\leq f\Big(\frac{x_1+\dots+x_n}{n}\Big).
}}

\begin{proof} Let $p\in\,]0,1]$. First we show that $F_W^p(S)$ is closed under the operation $\oplus$. To see this, let $f,g:S\to[0,1]$ be $W$-nondecreasing functions with $\sup f,\sup g\geq p$. If $x\in S$ and $w\in W$, then
\Eq{*}{
  (f\oplus g)(x+w)
  &=\sup_{\substack{u,v\in S\\ u+v=x+w}} \min\big(f(u),g(v)\big)
  \geq\sup_{\substack{u,v'\in S\\ u+v'=x}} \min\big(f(u),g(v'+w)\big)\\
  &\geq\sup_{\substack{u,v'\in S\\ u+v'=x}} \min\big(f(u),g(v')\big)=(f\oplus g)(x),
}
which shows that $f\oplus g$ is also $W$-nondecreasing. Let $\eta<p$ be arbitrary. Then there exist $u_0,v_0\in S$ such that $f(u_0)>\eta$ and $g(v_0)>\eta$ hold. Then we have  $(f\oplus g)(u_0+v_0)\geq\min\big(f(u_0),g(v_0)\big)>\eta$, proving that $\sup(f\oplus g)>\eta$. Taking te limit $\eta\to p$, this implies that $f\oplus g\in F_W^p(S)$.

The commutativity of the operation $\oplus$ is a consequence of the commutativity of the group operation $+$ in $G$. To verify the associativity, let $f,g,h\in F_W^p(S)$. Then, for all $x\in S$,
\Eq{*}{
  ((f\oplus g)\oplus h)(x)
  &=\sup_{\substack{u,v\in S\\ u+v=x}} \min\big((f\oplus g)(u),h(v)\big)
  =\sup_{\substack{u,v\in S\\ u+v=x}} \min\bigg(\sup_{\substack{s,t\in S\\ s+t=u}} \min\big(f(s),g(t)\big),h(v)\big)\bigg)\\
  &=\sup_{\substack{u,v\in S\\ u+v=x}} \sup_{\substack{s,t\in S\\ s+t=u}} \min\big(\min\big(f(s),g(t)\big),h(v)\big)\big)
  =\sup_{\substack{s,t,v\in S\\ s+t+v=x}}  \min\big(f(s),g(t),h(v)\big).
}
A similar argument shows that
\Eq{*}{
  (f\oplus (g\oplus h))(x)=\sup_{\substack{u,s,t\in S\\ u+s+t=x}}  \min\big(f(u),g(s),h(t)\big),
}
which results the desired equality $((f\oplus g)\oplus h)(x)=(f\oplus (g\oplus h))(x)$.

To see that the characteristic function $\chi_W$ of $W$ is a unital element of the semigroup $F_W^p(S)$, observe that $\chi_W$ is a $W$-nondecreasing function and, for all $x\in S$,
\Eq{*}{
  (f\oplus \chi_W)(x)
  =\sup_{\substack{u,v\in S\\ u+v=x}} \min\big(f(u),\chi_W(v)\big)
  =\sup_{\substack{u\in S,\,v\in W\\ u+v=x}} f(u)
  \leq\sup_{\substack{u\in S,\,v\in W\\ u+v=x}} f(u+v)=f(x).
}
On the other hand, by taking $v=0$, we can see that the inequality \Eq{*}{
f(x)\leq \sup_{\substack{u\in S,\,v\in W\\ u+v=x}} f(u)
}
holds, which finally implies the equality $(f\oplus \chi_W)(x)=f(x)$.

It is obvious that $(F_W^p(S),\leq)$ is a partially ordered set. We prove that the operation $\oplus$ is monotone with respect to the ordering $\leq$. Indeed, if $f,g,h\in F_W^p(S)$ and $g\leq h$ on $S$, then, for all $x\in S$,
\Eq{*}{
  (f\oplus g)(x)
  =\sup_{\substack{u,v\in S\\ u+v=x}} \min\big(f(u),g(v)\big)
  \leq\sup_{\substack{u,v\in S\\ u+v=x}} \min\big(f(u),h(v)\big)
  =(f\oplus h)(x).
}
So far we have shown that $(F_W^p(S),\oplus,\leq)$ is an ordered commutative semigroup.

In the rest of the proof, we prove that, this structure with the operation $\odot$ forms a cornet. 

First we show that $n\odot f\in F_W^p(S)$ whenever $n\in\N$ and $f\in F_W^p(S)$.
Indeed,
\Eq{*}{
  \sup_{x\in S} (n\odot f)(x)=\sup_{x\in S} f\Big(\frac{x}{n}\Big)
  \geq \sup_{y\in S} f\Big(\frac{n\cdot y}{n}\Big) = \sup_{y\in S}f(y)\geq p,
}
which proves that $\sup (n\odot f)\geq p$. On the other hand, if $x\leq_W y$, then $\frac{x}{n}\leq_W \frac{y}{n}$. By the $W$-nondecreasingness of $f$, this implies $f(\frac{x}{n})\leq f(\frac{y}{n})$, that is, $(n\odot f)(x)\leq n\odot f)(y)$. Therefore, $(n\odot f)$ is also $W$-nondecreasing. 

For property (i) of \defi{2}, let $n,m\in\N$ and $f\in F_W^p(S)$.
Then, for all $x\in S$,
\Eq{*}{
  ((nm)\odot f)(x)=f\Big(\frac{x}{nm}\Big)
  =f\Big(\frac{x/n}{m}\Big)=(m\odot f)\Big(\frac{x}{n}\Big)
  =(n\odot (m\odot f))(x),
}
which shows the expected identity $(nm)\odot f=n\odot (m\odot f)$. 

For property (ii) of \defi{2}, let $n\in\N$ and $f,g\in F_W^p(S)$. Then, for all $x\in S$,
\Eq{*}{
  (n\odot (f\oplus g))(x)
  &=(f\oplus g)\Big(\frac{x}{n}\Big)
  =\sup_{\substack{u,v\in S\\ u+v=\frac{x}{n}}} 
     \min\big(f(u),g(v)\big)
  =\sup_{\substack{u',v'\in S\\ u'+v'=x}} 
   \min\Big(f\Big(\frac{u'}{n}\Big),g\Big(\frac{v'}{n}\Big)\Big)\\
  &=\sup_{\substack{u',v'\in S\\ u'+v'=x}} 
     \min\big((n\odot f)(u'),(n\odot g)(v')\big)
  =((n\odot f)\oplus(n\odot g))(x),
}
which proves the equality $n\odot (f\oplus g)=(n\odot f)\oplus(n\odot g)$.

For property (iii) of \defi{2}, let $n,m\in\N$ and $f\in F_W^p(S)$.
Then, for all $x\in S$,
\Eq{*}{
  ((n+m)\odot f)(x)
  &=f\Big(\frac{x}{n+m}\Big)
  =\min\Big(f\Big(\frac{1}{n}\cdot\frac{nx}{n+m}\Big),
  f\Big(\frac{1}{m}\cdot\frac{mx}{n+m}\Big)\Big)\\
  &\leq\sup_{\substack{u,v\in S\\ u+v=x}} 
  \min\Big(f\Big(\frac{u}{n}\Big),f\Big(\frac{v}{m}\Big)\Big)
  =((n\odot f)\oplus(m\odot f))(x),
}
which shows the desired inequality $(n+m)\odot f\leq (n\odot f)\oplus(m\odot f)$.

For property (iv) of \defi{2}, let $n\in\N$ and $f,g\in F_W^p(S)$ with $f\leq g$. Then, for all $x\in S$,
\Eq{*}{
  (n\odot f)(x)=f\Big(\frac{x}{n}\Big)
  \leq g\Big(\frac{x}{n}\Big)=(n\odot g)(x),
}
which yields that $n\odot f\leq n\odot g$.

The property (v), which is the equality $1\odot f=f$, is obvious. The equality $n\odot\chi_W=\chi_W$ easily follows from the equivalence of the inclusions $\frac{x}{n}\in W$ and $x\in W$. Thus property (vi) of \defi{2} is also satisfied.

We now show that $(F_W^{p}(S),\leq)$ is a complete partially ordered set. Let $\F:=\{f_\gamma\mid\gamma\in\Gamma\}$ be a family of elements in $F_W^p(S)$ bounded from below by $g\in F_W^{p}(S)$. Define $f:S\to[0,1]$ by $f:=\inf_{\gamma\in\Gamma}f_\gamma$. We prove that $f$ is also a member of $F_W^p(S)$. By the inequality $g\leq f_\gamma$, it follows that $g\leq f$ and hence $p\leq \sup g\leq f$. Let $x,y\in S$ with $x\preceq_W y$, that is, with $y-x\in W$. Then, for all $\gamma\in\Gamma$ the $W$-nondecreasingness of $f_\gamma$ gives $f_\gamma(x)\leq f_\gamma(y)$. Taking the infimum with respect to $\gamma\in\Gamma$ side by side, it follows that $f(x)\leq f(y)$, which proves that $f$ is also $W$-nondecreasing and hence $f\in F_W^p(S)$. Clearly, $f$ is the infimum of the family $\F$ and this shows that $(F_W^p(S),\leq)$ is a complete partially ordered set.

We verify that the map $\Phi$ defined by \eq{Phi} is an injective cornet-preserving mapping of $(P_W(S),+,*,\subseteq)$ into $(F_W^1(S),\oplus,\odot,\leq)$. Clearly, if $A\in P_W(S)$, then $\Phi(a)=\chi_A$ is $W$-nondecreasing and $\sup\Phi(A)=\sup\chi_A=1$, which shows that $\Phi(A)\in F_W^1(S)$. The injectivity of $\Phi$ is obvious. To prove that $\Phi$ preserves the addition, let $A,B\in P_W(S)$. Then, for $x\in S$, it is easy to see that 
\Eq{*}{
  \sup_{\substack{u,v\in S\\ u+v=x}}\min(\chi_A(u),\chi_B(v))=1
}
if and only if there exist $u\in A$, $v\in B$ such that $x=u+v$, that is, if $x\in A+B$. This proves that, for all $x\in S$, 
\Eq{*}{
  (\chi_A\oplus\chi_B)(x)=\sup_{\substack{u,v\in S\\ u+v=x}}\min(\chi_A(u),\chi_B(v))
  =\chi_{A+B}(x).
}
As a consequence of this equality, it follows that $\Phi(A+B)=\Phi(A)\oplus\Phi(B)$.

Let $A,B\in P_W(S)$ and $n\in\N$. It is clear that
\Eq{nA}{
  \{n\cdot a\mid a\in A\}
  \subseteq n*A=\{n\cdot a+w\mid a\in A,\,w\in W\}.
}
In fact, this inclusion is an equality. To see this, let $x\in S$
be of the form $x=n\cdot a+w$ for some $a\in A$ and $w\in W$.
Then, by the divisibility of $W$, we have that $w/n\in W$. Thus, the $W$-invariance of $A$ yields that $a'=a+(w/n)\in A$ and hence
$x$ is of the form $n\cdot a'$ for some element $a'\in A$, which shows that it belongs to the left hand side set in \eq{nA}.

Using the equality \eq{nA}, for $x\in S$, we have
\Eq{*}{
  \chi_{n*A}(x)=\chi_{A}\Big(\frac{x}{n}\Big)=n\odot\chi_A(x),
}
which proves the equality $\Phi(n*A)=n\odot\Phi(A)$. 

If $A,B\in P_W(S)$ with $A\subseteq B$, then $\chi_A(x)\leq\chi_B(x)$ holds for all $x\in S$, which shows that $\Phi(A)\leq\Phi(B)$, that is, $\Phi$ preserves the ordering as well.

To prove the last assertion of the proposition, assume that $f\in F_W(S)$ is an $n$-convex element. Let $x_1,\dots,x_n\in S$. By the $n$-convexity of $f$, we have that $n\cdot f\leq n\odot f$, that is, for all $x\in S$,
\Eq{*}{
 \sup_{\substack{u_1,\dots,u_n\in S\\ u_1+\dots+u_n=x}} 
    \min\big(f(u_1),\dots,f(u_n)\big)
    =(f\oplus\dots\oplus f)(x)
 \leq (n\odot f)(x)=f\Big(\frac{x}{n}\Big).
}
By taking $x:=x_1+\dots+x_n$, with $u_1:=x_1,\dots,u_n:=x_n$, it follows that \eq{qc} holds. The proof of the reversed implication is analogous.
\end{proof}

\section{Topological notions and boundedness in cornets}

In a natural way, we can introduce the notions of nonnegative and Archimedean elements in a cornet with the following definition.

\Defi{4}{
In a cornet $(X,+,*,\preceq)$ an element $x\in X$ is said to be \emph{nonnegative} if $0\preceq x$ holds. The element $x$ is called \emph{Archimedean}, denoted by $0\prec x$, if, for all $u\in X$, there exists $n_0\in\N$ such that $0\preceq u+n*x$ for all $n_0\leq n$. The set of all nonnegative and Archimedean elements in $X$ will be denoted by $X_\preceq$ and $X_\prec$, respectively.
}

The properties of nonnegative and Archimedean elements are established in the following assertion.

\Prp{1}{
Let $(X,+,*,\preceq)$ be a cornet. Then $X_\prec$ is contained in $X_\preceq$ and
\Eq{incl}{
  X_\prec+X_\preceq\subseteq X_\prec.
}
In addition, $X_\prec$ and $X_\preceq$ are subcornets of $(X,+,*,\preceq)$.}

\begin{proof} To show that $X_\prec\subseteq X_\preceq$, let $x\in X_\prec$. Then, there exists $n_0\in\N$ such that $0\preceq n_0*x$. Therefore, $n_0*0\preceq n_0*x$, which implies that $0\preceq x$, i.e., $x\in X_\preceq$.

Let $x\in X_\prec$ and $y\in X_\preceq$. Then, for any $u\in X$, there exists $n_0\in\N$ such that $0\preceq u+n*x$ for all $n_0\leq n$. On the other hand, $0=n*0\preceq n*y$, therefore $0\preceq u+n*(x+y)$ holds for all $n_0\leq n$, which implies that $x+y$ is Archimedean and proves the inclusion \eq{incl}. 

If $x$ and $y$ are nonnegative elements, then by the axioms of cornets, $0\preceq y=0+y\preceq x+y$, which shows that $X_\preceq$ is closed under addition. Similarly, if $x\in X_\preceq$ and $m\in\N$, then $0=m*0\preceq m*x$ proving that $X_\preceq$ is closed under $*$-multiplication and hence $(X_\preceq,+,*,\preceq)$ is a subcornet of $(X,+,*,\preceq)$.

By \eq{incl} we have that $X_\prec+X_\prec\subseteq X_\prec+X_\preceq\subseteq X_\prec$, which shows that $X_\prec$ is closed under addition.
To prove that $X_\prec$ is also closed under $*$-multiplication, let $x\in X_\prec$ and $m\in\N$ and let $u\in X$ be arbitrary. Then there exists $n_0\in\N$ such that $0\preceq u+n*x$ for all $n\geq n_0$.
In particular, we have that $0\preceq u+(km)*x=u+k*(m*x)$ for all $k\geq \frac{n_0}{m}$, which shows that $m*x$ is also Archimedean. Hence, $(X_\prec,+,*,\preceq)$ is also a subcornet of $(X,+,*,\preceq)$
\end{proof}

In what follows, we introduce the notions of continuity of the addition, boundedness and closedness with respect to a subsemigroup of Archimedean elements. For comparison, we recall first the standard topological concepts for abelian groups.

\Defi{5-}{If $(G,+)$ is an abelian group and $\T$ is a Hausdorff topology on $G$, then we say that $(G,\T,+)$ is a \emph{topological group} if the $(x,y)\mapsto x-y$ is a continuous map of $G\times G$ into $G$. A subset $U\subseteq G$ is said to be \emph{convex} if, for all $n\in\N$, 
\Eq{*}{
\{u_1+\dots+u_n\mid u_1,\dots,u_n\in U\}=\{n\cdot u\mid u\in U\}.
} 
We say that $G$ is \emph{locally convex} if every neighborhood of $0$ contains a convex neighborhood of $0$. A subset $H\subseteq G$ is said to be \emph{topologically bounded} if, for all neighborhood $U$ of $0$, there exists $n\in\N$ such that $H\subseteq \{u_1+\dots+u_n\mid u_1,\dots,u_n\in U\}$.}

%Let $(X,+,*,\preceq)$ be a cornet and let $\A$ be a subsemigroup of $X_\prec$.
%We say that 

\Defi{6-}{Let $(X,+,*,\preceq)$ be a cornet and let $\A$ be a subsemigroup of $X_\prec$. We say that the \emph{addition is $\A$-continuous} if, for all $a\in\A$, there exists $b\in\A$ such that $b+b\preceq a$ holds.}

\Lem{Atop}{Let $(X,+,*,\preceq)$ be a cornet, let $\A$ be a subsemigroup of $X_\prec$ and assume that the addition is $\A$-continuous. Then the $\cdot$-multiplication and the $*$-multiplication are $\A$-continuous, that is, for all $a\in\A$ and $n\in\N$, there exists $b\in\A$ such that $n\cdot b\preceq a$ and $n*b\preceq a$, respectively.}

\begin{proof}
 Let $a\in\A$. By the continuity of addition there exists $a_1\in\A$ such that $a_1+a_1\preceq a$, or in other words $2\cdot a_1\preceq a$. Applying the continuity of addition again for $a_1\in\A$, we get that there exists $a_2\in\A$ such that $2\cdot a_2\preceq a_1$, which implies that $4\cdot a_2\preceq a$. Continuing this process, we can construct a sequence $a_k\in\A$ such that $2^k\cdot a_k\preceq a$ holds for all $k\in\N$. Now let $n\in\N$ and choose $k\in\N$ so that $n\leq 2^k$. Then $n\cdot a_k=n\cdot a_k+(2^k-n)\cdot0\preceq 2^k\cdot a_k\preceq a$. Thus, $n\cdot b\preceq a$ holds with $b:=a_k$. The inequality $n*b\preceq n\cdot b$ implies that $n*b\preceq a$ is also valid.
\end{proof}

\Defi{5}{Let $(X,+,*,\preceq)$ be a cornet and let $\A$ be a subsemigroup of $X_\prec$. We say that an element $x\in X$ is \emph{$\A$-bounded} if, for all $a\in\A$, there exists $n_0\in\N$ such that $x\preceq n*a$ for all $n_0\leq n$.}	

\Prp{atop}{
        Let $(X,+,*,\preceq)$ be a cornet and let $\A$ be a subsemigroup of $X_\prec$ such that the addition is $\A$-continuous. Then the $\A$-bounded elements form a subcornet of $(X,+,*,\preceq)$. 
}

\begin{proof}
	Let $x,y$ be $\A$-bounded elements of $X$ and let $a\in\A$ be arbitrary. Since $X$ is $\A$-topological therefore there exists $b\in\A$ such that $b+b\preceq a$. Using the $\A$-boundedness property of $x,y$, there exist $k_0,m_0\in\N$ such that
	\Eq{*}{
		x\preceq k*b\qquad &\text{if $k\geq k_0$},\\
		y\preceq m*b\qquad &\text{if $m\geq m_0$.}
	}
	Thus,
	\Eq{*}{
		x\preceq n*b\qquad &\text{if $n\geq \max(k_0,m_0)$},\\
		y\preceq n*b\qquad &\text{if $n\geq \max(k_0,m_0)$.}
	}
	Adding the above inequalities side by side, we get
	\Eq{*}{
		x+y\preceq n*(b+b)\preceq n*a\qquad \text{if $n\geq \max(k_0,m_0)$,}
	}
	which proves that the set of $\A$-bounded elements is closed under addition.
	
	To prove that the set of $\A$-bounded elements is closed under the $*$-multiplication, let $x\in X$ be $\A$-bounded, $m\in\N$ and $a\in\A$. Then, there exist $b\in\A$ such that $m*b\preceq a$. Then there exists $n_0\in\N$ such that $x\preceq n*b$ holds for all $n\geq n_0$. Using the monotonicity property of $*$-multiplication, for all $n\geq n_0$, it follows that $m*x\preceq m*(n*b)=(mn)*n=n*(m*b)\preceq n*a$. This shows that the set of $\A$-bounded elements is closed under the $*$-multiplication and proves that $\A$-bounded elements form a subcornet of $(X,+,*,\preceq)$.
\end{proof}

\Defi{6+}{
	Let $(X,+,*,\preceq)$ be a cornet and let $\A$ be a subsemigroup of $X_\prec$. Given an element $x\in X$, we say that $y\in X$ is the \emph{$\A$-closure} of $x$ if $y\preceq x+a$ holds for all $a\in A$ and, $y$ is the largest element of $X$ with is property, i.e., if $z\preceq x+a$ holds for all $a\in A$, then $z\preceq y$.
	It is clear that the $\A$-closure of an element, if exists, is unique and is denoted by $\cl_\A(x)$. An element $x$ is called \emph{$\A$-closed} if $x=\cl_\A(x)$. The set of all elements of $X$ which possess an $\A$-closure will be denoted by $\Cl_\A$.
}

\Prp{cl}{Let $(X,+,*,\preceq)$ be a cornet and let $\A$ be a subsemigroup of $X_\prec$ such that $X$ is $\A$-topological. Then we have the following assertions.
\begin{enumerate}[(i)]
 \item If $x\in\Cl_\A$, then $x\preceq\cl_\A(x)$.
 \item If $x,y\in \Cl_\A$ and $x\preceq y$, then $\cl_\A(x)\preceq \cl_\A(y)$.
 \item If $x\in \Cl_\A$, then $\cl_\A(x)\in\Cl_\A$ and $\cl_\A(x)=\cl_\A(\cl_\A(x))$.
 \item If $x\in\Cl_\A$ and $x\preceq y\preceq \cl_A(x)$, then $y\in\Cl_\A$ and $\cl_\A(y)=\cl_\A(x)$.
 \item If $x,y\in \Cl_\A$ with $x+y\in \Cl_\A$, then $\cl_\A(x)+\cl_\A(y)\in\Cl_\A$ and $\cl_\A(\cl_\A(x)+\cl_\A(y))=\cl_\A(x+y)$.
 \item If $x\in \Cl_\A$, $n\in\N$ and $n\cdot x\in\Cl_\A$, then $n\cdot \cl_\A(x)\in\Cl_\A$ and $\cl_\A(n\cdot\cl_\A(x))=\cl_\A(n\cdot x)$.
 \item If $x\in \Cl_\A$, $n\in\N$ and $n*x\in\Cl_\A$, then $n*\cl_\A(x)\in\Cl_\A$ and $\cl_\A(n*\cl_\A(x))=\cl_\A(n*x)$.
 \item If $x\in \Cl_\A$ is $\A$-bounded, then $\cl_\A(x)$ is also $\A$-bounded.
 \item If $n\in\N$, $x\in \Cl_\A$ is $n$-convex, $n*x\in\Cl_\A$ and $n*\cl_\A(x)$ is $\A$-closed, then $\cl_\A(x)$ is also $n$-convex.
\end{enumerate}
}

\begin{proof} \hfill\null

(i) Let $x\in \Cl_\A$. Clearly $x\preceq x+a$ holds for all $a\in\A$. This implies that $x\preceq \cl_\A(x)$.
	
(ii) Let $x,y\in \Cl_\A$ such that $x\preceq y$. By using the definition of $\A$-closure, for all $a\in\A$, we have that
	\Eq{*}{
		\cl_\A(x)\preceq x+a\preceq y+a,
	}
	which implies that $\cl_\A(x)\preceq \cl_\A(y)$.

(iii) Let $x\in\Cl_\A$. We need to show that $y=\cl_A(x)$ is $\A$-closed. The inequality $y\preceq y+a$ is trivial for all $a\in\A$. Assume that $z\preceq y+a$ holds for all $a\in\A$. Let $b\in\A$ be arbitrary. By the $\A$-continuity of the addition, there exist $a,c\in\A$ such that $a+c\preceq b$. We have that $y\preceq x+c$, hence, 
\Eq{*}{
z\preceq y+a\preceq (x+c)+a\preceq x+b.
}
Since $b$ was arbitrary, this implies that $z\preceq \cl_\A(x)=y$, which shows that $y$ is the $\A$-closure of itself.

(iv) Let $x\in\Cl_\A$ and $x\preceq y\preceq \cl_A(x)$. We need to show that $\cl_\A(x)$ is the $\A$-closure of $y$. 

On one hand, for all $a\in\A$, we have $\cl_\A(x)\preceq x+a\preceq y+a$. On the other hand, assume that $z\in X$ satisfies $z\preceq y+a$ for all $a\in\A$. Then $z\preceq \cl_\A(x)+a$ for all $a\in\A$, which, by property (iii) implies that $z\preceq\cl_A(x)$. This shows that $\cl_\A(x)$ is the $\A$-closure of $y$. 
	
(v) Let $x,y\in X$ be arbitrary. First we show that $\cl_\A(x)+\cl_\A(y)\preceq x+y+a$ holds for all $a\in\A$. Indeed, if $a\in\A$, then there exist $b,c\in\A$ such that $b+c\preceq a$ and we get
\Eq{*}{
	\cl_\A(x)+\cl_\A(y)\preceq (x+b)+(y+c)\preceq x+y+a.
	}
This inequality and property (i) imply that $x+y\preceq\cl_\A(x)+\cl_\A(y)\preceq \cl_A(x+y)$. Applying properties (iii) and (iv), it follows that $\cl_\A(x)+\cl_\A(y)\in\Cl_\A$ and $\cl_A(\cl_\A(x)+\cl_\A(y))=\cl_A(x+y)$, which was to be proved.
	
(vi) Let $x\in \Cl_\A$, $n\in\N$ with $n\cdot x\in\Cl_\A$. First we show that $n\cdot\cl_\A(x)\preceq n\cdot x+a$ holds for all $a\in\A$. Indeed, if $a\in\A$, then there exist $b\in\A$ such that $n\cdot b\preceq a$. Therefore, the inequality $\cl_\A(x)\preceq x+b$ implies that
\Eq{*}{
 n\cdot\cl_\A(x)\preceq n\cdot(x+b)=n\cdot x+n\cdot b\preceq n\cdot x+a.
}
In view of property (i) and this inequality, we get that $n\cdot  x\preceq n\cdot\cl_\A(x)\preceq \cl_A(n\cdot x)$. Applying properties (iii) and (iv), it follows that $n\cdot\cl_\A(x)\in\Cl_\A$ and $\cl_A(n\cdot\cl_\A(x))=\cl_A(n\cdot x)$, which yields the statement.
	
(vii) Let $x\in \Cl_\A$, $n\in\N$ with $n*x\in\Cl_\A$. First we show that $n*\cl_\A(x)\preceq n*x+a$ holds for all $a\in\A$. Indeed, if $a\in\A$, then there exist $b\in\A$ such that $n*b\preceq a$. Therefore, the inequality $\cl_\A(x)\preceq x+b$ implies that
\Eq{*}{
 n*\cl_\A(x)\preceq n*(x+b)=n*x+n*b\preceq n*x+a.
}
In view of property (i) and this inequality, we get that $n*x\preceq n*\cl_\A(x)\preceq \cl_A(n*x)$. Applying properties (iii) and (iv), it follows that $n*\cl_\A(x)\in\Cl_\A$ and $\cl_A(n*\cl_\A(x))=\cl_A(n*x)$, which yields the statement.
	
(viii) Let $x\in X$ be an $\A$-bounded element and let $a\in\A$ be arbitrary. Then, for some $b,c\in\A$, we have that $b+c\preceq a$. By the $\A$-boundedness of $x$, there exists $n_0\in\N$ such that $x\preceq n*b$ holds for all $n\geq n_0$. On the other hand, $c\in X_\prec$, hence $n*c\in X_\prec$ for all $n\in\N$. Consequently, for all $n\geq n_0$, we have that
\Eq{*}{
  \cl_\A(x)\preceq x+n*c\preceq n*b+n*c=n*(b+c)\preceq n*a,
}
which shows that $\cl_\A(x)$ is also an $\A$-bounded element of $X$.
	
(ix) Let $x\in X$ be $n$-convex such that $n*x\in\Cl_\A$ and $n*\cl_\A(x)$ is $\A$-closed. Then, the equality $n\cdot x=n*x$ yields that $n\cdot x\in\Cl_\A$ and properties (i), (ii), (v) and (vi), (vii) imply that
	\Eq{*}{
	  n\cdot\cl_\A(x)\preceq\cl_\A(n\cdot\cl_\A(x))=\cl_\A(n\cdot x)=\cl_\A(n*x)=\cl_\A(n*\cl_\A(x)).
	}
	Using that $n*\cl_\A(x)$ is $\A$-closed, this inequality implies $n\cdot\cl_\A(x)\preceq n*\cl_\A(x)$. The reversed inequality holds automatically, hence $n\cdot\cl_\A(x)= n*\cl_\A(x)$, which proves the $n$-convexity of $\cl_\A(x)$.
\end{proof}

In what follows, we investigate the connection among the notions of boundedness, closedness and convexity. 

In the subsequent propositions, we consider the cornets that were introduced in Proposition \ref{P1+}, \ref{P2}, \ref{P3} and we determine the Archimedean, the bounded and closed elements in these structure. 

\Prp{3.1}{Let $(G,+)$ be a topological abelian group such that there is no proper open subgroup of $G$. Let $W\subseteq G$ be a wedge with $W^\circ\neq\emptyset$ and let $S$ be a subsemigroup of $G$ containing $W$. Then we have the following claims:
\begin{enumerate}[(i)]
	\item In the cornet $(S,+,\cdot,\preceq_W)$ the set of nonnegative elements equals $W$.
	\item The set $W^\circ$ is a subsemigroup of the Archimedean elements.
	\item Every element of $S$ is $W^\circ$-bounded.
	\item If, in addition, $G$ is locally convex and $W$ is topologically closed and $W^\circ=W^\circ+W^\circ$, then every element of $S$ is also $W^\circ$-closed.
\end{enumerate}
}

\begin{proof} (i) The nonnegativity of an element $x\in S$ with respect to the ordering $\preceq_W$, by definition, means that $x=x-0\in W$. This proves that $W$ equals the set of nonnegative elements.

(ii) We have that $W^\circ+W^\circ\subseteq W+W\subseteq W$ and $W^\circ+W^\circ$ is also open, therefore, $W^\circ+W^\circ\subseteq W^\circ$ showing that $W^\circ$ is a subsemigroup of $G$.

To show that every element of $W^\circ$ is Archimedean, let $x\in W^\circ$ be arbitrary and define
\Eq{Tx}{
  T_x:=\bigcup_{n\in\N} (W^\circ-n\cdot x)\cap(n\cdot x-W^\circ).
}
We show that $T_x$ is an open subgroup of $G$. The opennes of $T_x$ is obvious. Let $y,z\in T_x$. Then there exist $n,m\in\N$ and $v_1,v_2,w_1,w_2\in W^\circ$ such that
\Eq{*}{
   y=v_1-n\cdot x=n\cdot x-v_2,\qquad z=w_1-m\cdot x=m\cdot x-w_2.
}
Therefore,
\Eq{*}{
y-z&=(v_1+w_2)-(n+m)\cdot x=(n+m)\cdot x-(v_2+w_1),
}
and hence
\Eq{*}{
  y-z\in (W^\circ-(n+m)\cdot x)\cap((n+m)\cdot x-W^\circ)\subseteq T_x,
}
which proves that $T_x$ is an open subgroup of $G$. Thus, $T_x$ cannot be proper, in other words, $T_x=G$. 

Let $u\in S$ be arbitrary. Then $S\subseteq T_x$ implies that there exists
$n_0\in\N$ such that $u\in W^\circ -n_0\cdot x$. If $n>n_0$, then $(n-n_0)\cdot x\in (n-n_0)\cdot W^\circ\subseteq W^\circ$. Therefore,
\Eq{*}{
   u+n\cdot x=u+n_0\cdot x+(n-n_0)\cdot x\in (W^\circ-n_0\cdot x)+n_0\cdot x+W^\circ\subseteq W^\circ\subseteq W,
}
i.e., $0\preceq_W u+n\cdot x$, which proves that any element of $W^\circ$ is Archimedean.

(iii) Now we are going to show that every element of $S$ is $W^\circ$-bounded.
Let $u\in S$ be fixed and $x\in W^\circ$ be arbitrary. As we have seen above, the set $T_x$ defined by \eq{Tx}, covers $S$, therefore, there exists $n_0\in\N$ such that $u\in n_0\cdot x-W^\circ$. If $n>n_0$, then $(n-n_0)\cdot x\in (n-n_0)\cdot W^\circ\subseteq W^\circ$. Thus,
\Eq{*}{
  n\cdot x-u=(n-n_0)\cdot x+n_0\cdot x-u\in W^\circ+W^\circ\subseteq W^\circ\subseteq W,
}
which shows that $u\preceq_W n\cdot x$ holds for all $n>n_0$. This proves that $u$ is $W^\circ$-bounded.

(iv) Finally, assume that $W$ is topologically closed. We are going to verify that every element $u$ of $S$ is $W^\circ$-closed, that is, $u$ is the greatest lower bound of the set $\{u+x\mid x\in W^\circ\}$. It is clear that $u$ is a lower bound. Assume that $v\in S$ is a lower bound for $\{u+x\mid x\in W^\circ\}$, that is, $u-v+x\in W$ holds for all $x\in W^\circ$. Let $x\in W^\circ$ be fixed and $n\in\N$. 

Assume that $W^\circ=W^\circ+W^\circ$. Then 
\Eq{*}{
  W^\circ=\{x_1+\dots+x_n\mid x_1,\dots,x_n\in W^\circ\}.
}
Therefore, there exist $x_1,\dots,x_n\in W^\circ$ such that $x=x_1+\dots+x_n$.
This implies that 
\Eq{*}{
n\cdot(u-v)+x=(u-v+x_1)+\cdots+(u-v+x_n)\in W.
}

If $u-v\not\in W$, then $u-v\in G\setminus W$. Then there exists an open convex and symmetric neighborhood of $0$ such that $u-v+U\subseteq G\setminus W$. Consider the set $S$ defined by
\Eq{*}{
  S:=\bigcup_{n=1}^\infty\{y_1+\dots+y_n\mid y_1,\dots,y_n\in U\}.
}
Then $S$ is an open subgroup of $G$, hence $S=G$, which yields that $x\in S$. This implies that, for some $n\in\N$ and $y_1,\dots,y_n\in U$ the equality $x=y_1+\dots+y_n$ holds. By the convexity of $U$, it follows that there exists $y\in U$, such that $x=n\cdot y$. Thus
\Eq{*}{
  n\cdot(u-v+y)=n\cdot(u-v)+x\in W,
}
which yields that $u-v+y\in W$ contradicting $u-v+y\in u-v+U \subseteq G\setminus W$. This contradiction shows that $u-v\in W$ must be valid, i.e., $v\preceq_W u$ holds.
\end{proof}

\Prp{3.2}{
	Let $(G,+)$ be a topological abelian group such that there is no proper open subgroup of $G$. Let $W$ be a wedge and let $S$ be a subsemigroup of $G$ containing $W$. Let $P_W(S)$ denote the collection of all nonempty $W$-invariant subsets of $S$. Define the operations $+$ and $*$ by \eq{+*}. Then the following statements hold:
	\begin{enumerate}[(i)]
		\item The set of nonnegative elements of the cornet $(P_W(S),+,*,\subseteq)$ consists of those $W$-invariant subsets of $S$ that contain $0$ (which denotes the neutral element of $G$).
		\item The collection $\A$ of those $W$-invariant subsets which contain an open convex neighborhood $C\in P_W(S)$ of $0$ is a subsemigroup of the Archimedean elements.
		\item An element of $P_W(S)$ is $\A$-bounded if it is the sum of a topologically bounded subset of $S$ and $W$.
		\item If, in addition, $G$ is locally convex, then any topologically closed element of $P_W(S)$ is also $\A$-closed and the addition is $\A$-continuous.
	\end{enumerate}
}

\begin{proof} (i) The unit element of $(P_W(S),+,*,\subseteq)$ is the set $W$.
Now, an element $A\in P_W(S)$ is nonnegative if $W\subseteq A$. This inclusion is equivalent to the condition $0\in A$ because $A$ is $W$-invariant. 

(ii) Let $A,B\in\A$. Then there exist open convex sets $C,D\in P_W(S)$ such that $0\in C\subseteq A$ and $0\in D\subseteq B$, Then, by \lem{C}, the set $C+D$ is convex and also open, furthermore $0\in C+D\subseteq A+B$. Therefore, $A+B\in\A$, which proves that $(\A,+)$ is a semigroup.

To prove that the elements of $\A$ are Archimedean, let $A\in\A$. Then there exists an open and convex $C\in P_W(S)$ such that $0\in C\subseteq A$. Define
\Eq{TA}{
  T_C:=\bigcup_{n\in\N} (n*C)\cap(-n*C).
}
We prove that $T_C$ is a subgroup of $G$. Let $x,y\in T_C$. Then there there exist $n,m\in\N$ such that $x\in (n*C)\cap(-n*C)$ and $y\in(m*C)\cap(-m*C)$. Using the definition of the $*$ multiplication, it follows that there exist $a_1,a_2,b_1,b_2\in C$ and $v_1,v_2,w_1,w_2\in W$ such that
\Eq{*}{
  x=n\cdot a_1+v_1=-n\cdot a_2-v_2, \qquad
  y=m\cdot b_1+w_1=-m\cdot b_2-w_2.
}
By the $(n+m)$-convexity of $C$, it follows that here exist $c_1,c_2\in C$
and $u_1,u_2\in W$ such that 
\Eq{*}{
  n\cdot a_1+m\cdot b_2=(n+m)\cdot c_1+u_1,\qquad
  n\cdot a_2+m\cdot b_1=(n+m)\cdot c_2+u_2.
}
Then 
\Eq{*}{
  x-y&=(n\cdot a_1+m\cdot b_2)+(v_1+w_2)
     \subseteq (n+m)\cdot c_1 +u_1+W\subseteq (n+m)*C,\\
  y-x&=(n\cdot a_2+m\cdot b_1)+(v_2+w_1)
     \subseteq (n+m)\cdot c_2 +u_2+W\subseteq (n+m)*C.
}
Therefore, $x-y\in T_C$, which proves that $T_C$ is a subgroup of $G$. The openness of $C$ implies that $n\cdot C$ is open for every $n\in\N$. Thus, by the convexity of $C$, we have that $n*C$ is open for every $n\in\N$. Consequently, $T_C$ is open and hence, by our assumption, $T_C$ is equal to $G$.

Before proving the further assertions, we show that $n*C\subseteq m*C$ if $n\leq m$. Indeed, let $u\in n*C$. Then, for some $c\in C$ and $w\in W$, we have $u=n\cdot c+w$. If $n\leq m$, then, the $m$-convexity of $C$ and $0\in C$ yield that there exist $d\in C$ and $v\in W$ such that $n\cdot c+(m-n)\cdot 0=n\cdot d+v$. Thus
\Eq{*}{
  u=n\cdot c+(m-n)\cdot 0+w=m\cdot d+v+w=m*C,
}
which verifies $n*C\subseteq m*C$.

Let $U\in P_W(S)$ be arbitrary and choose a fixed element $u\in U$. The inclusion $u\in T_C$ yields, for some $n_0\in\N$, that $-u\in n_0*C$. For $n\geq n_0$, this implies that $-u\in n*C\subseteq n*A$. Hence $0\in U+n*A$ holds for all $n\geq n_0$ which, according to the the first assertion, means that $U+n*A$ is nonnegative for all $n\geq n_0$. This proves that $A$ is Archimedean.

(iii) Let $B\in P_W(S)$ be the sum of a topologically bounded set $D\subseteq S$ and $W$. Let $A\in\A$ be fixed. Then there exists an open convex set $C\in P_W(S)$ such that $0\in C\subseteq A$. 

By the topological boundedness of $D$, we can find a number $n_0$ such that $D\subseteq n_0\cdot C$. Since $0\in C$, this implies that $D\subseteq n\cdot C$ for all $n\geq n_0$. By the convexity of $C$, this yields that $D\subseteq n*C$.  Consequently
\Eq{*}{
  B=D+W\subseteq n*C+W=n*C\subseteq n*A,
}
which proves that $B$ is $\A$-bounded.

%Now assume that $B$ is $\A$-bounded. Define
%\Eq{*}{
%  D:=\{b\in B\mid (b-W)\cap B=\{b\}\}.
%}
%$B=D+W$ ? $D$ is topologically bounded ?

(iv) In this part of the proof, we assume that $G$ is a locally convex topological group. Let $D\in P_W(S)$ be a topologically closed set. We need to show that $D$ is the $\A$-closure of itself. The inclusion $D\subseteq D+A$ trivially holds for all $A\in\A$ because $0\in A$. Assume now that, for some $E\in P_W(S)$, the inclusion $E\subseteq D+A$ holds for all $A\in\A$. We need to show that $E\subseteq D$. 

Let $e\in E$ be arbitray and assume that $e\not\in D$, i.e., $0\not\in e-D$. Using that $D$ is topologically closed, we have that $e-D$ is closed. Hence $0$ is an interior point of its complement. Thus there exists an open convex neighborhood $C_0$ of zero such that $C_0\cap(e-D)=\emptyset$. We show that then $(C_0+W)\cap(e-D)=\emptyset$. Indeed, if this not true, then
there exist $c\in C_0,w\in W$ and $d\in D$ such that $c+w=e-d$. Then, by the $W$-invariance of $D$, we get $c=e-(d+w)\in e-D$ which contradicts $C_0\cap(e-D)=\emptyset$. Then $A:=C_0+W$ is an open convex $W$-invariant neighborhood of zero which is disjoint from $e-D$. This implies that $e\not\in A+D$ which contradicts that $E\subseteq D+A$ holds for all $A\in\A$

Finally, we prove that the addition is $\A$-continuous. Let $A\in\A$. We need to show that there exists $B\in\A$ such that $B+B\subseteq A$. 

By $A\in\A$, there exists a convex neighborhood $C\in P_W(S)$ of $0$ such that $C\subseteq A$. By the continuity of the addition in $G$, there exists a neighborhood $D$ of $0$ such that $D+D\subseteq C$. Using the local convexity of the topology of $G$, we may assume that $D$ convex. Define $B$ as $D+W$. Then $B$ is a $W$-invariant convex neighborhood of $0$, hence $B\in\A$, and $B+B=(D+W)+(D+W)\subseteq C+W\subseteq C\subseteq A$. 
\end{proof}

\Prp{3.3}{
	Let $(G,+)$ be a topological abelian group such that there is no proper open subgroup of $G$. Let $W$ be a wedge and let $S$ be a uniquely divisible subsemigroup of $G$ containing $W$. Let, for $p\in\,]0,1]$, the set $F^p_W(S)$ be defined by \eq{Fp} and define the operations $\oplus$ and $\odot$ by \eq{oper}. Then, the following statements hold.
	\begin{enumerate}[(i)]
		\item The set of nonnegative elements of the cornets $(F^p_W(S),\oplus,\odot,\leq)$ consists of those $W$-invariant functions $f$ such that $f(0)=1$ (here $0$ denotes the neutral element of $G$).
		\item The cornet $(F^p_W(S),\oplus,\odot,\leq)$ has no Archimedean elements for $p\in\,]0,1[\,$. On the other hand, the collection $\A$ of those $a\in F^1_W(S)$ or which there exists an open convex neighborhood $C$ of $0$ such that $a|_C=1$ is a subsemigroup of the Archimedean elements of $(F^1_W(S),\oplus,\odot,\leq)$.
		\item An element $f$ of $F^1_W(S)$ is $\A$-bounded if\, $\supp(f):=\{u\in S\mid f(u)>0\}$ is covered by the sum of a topologically bounded subset of $S$ and $W$.
		\item If, in addition, $G$ is locally convex, then any upper semicontinuous element of $F^1_W(S)$ is also $\A$-closed and the addition $\oplus$ is $\A$-continuous.
	\end{enumerate}
}

\begin{proof}
(i) Let $p\in\,]0,1]$. The unit element of $(F^p_W(S),\oplus,\odot,\leq)$ is the characteristic function of $W$. Therefore, by definition, an element $f\in F^p_W(S)$ is nonnegative if $\chi_W\leq f$ which implies that $1=\chi_W(0)\leq f(0)$, whence $f(0)=1$ follows. On the other hand, if $f(0)=1$, then $W$-nondecreasingness implies of $f$ implies that $1=f(0)\leq f(w)$. Thus $\chi_W\leq f$ holds.
 
(ii) Let $p\in\,]0,1[\,$. Assume that $g\in F^p_W(S)$ is an Archimedean element of $F^p_W(S)$. Let $f:S\to[0,1]$ be a constant function with a value $p\leq f(0)<1$. Then $f\in F^p_W(S)$, on the other hand, for all $n\in\N$,
\Eq{*}{
  (f\oplus (n\odot g))(0)
  =\sup_{\substack{u,v\in S\\ u+v=0}}\min\big(f(u),(n\odot g)(v)\big)\leq f(0)<1
}
According to the first assertion of this theorem, $f\oplus (n\odot g)$ is not nonnegative for all $n\in\N$, which shows that $g$ cannot be Archimedean.

Let $f,g\in\A$. Then there exist open convex neighborhoods $C,D$ of $0$ such that $f|_C=1$ and $g|_D=1$. Then, for $x\in C+D$, we get
\Eq{*}{
  (f\oplus g)(x)
  =\sup_{\substack{u,v\in S\\ u+v=x}}\min\big(f(u),g(v)\big)
  \geq \sup_{\substack{u\in C,v\in D\\ u+v=x}}\min\big(f(u),g(v)\big)
  =1.
}
On the other hand, $C+D$ is an open convex neighborhood of $0$ on which $f\oplus g=1$. Therefore, $f\oplus g\in\A$.

Now we show that the elements of $\A$ are Archimedean in $(F^1_W(S),\oplus,\odot,\leq)$. Let $f\in\A$ and $g\in F^1_W(S)$ be arbitrary. We need to show that there exists $n_0\in\N$ such that $(g\oplus(n\odot f))(0)=1$ holds for all $n\in\N$. 
Let $0<\eta<1$ be arbitrary and choose $u_0\in S$ such that $g(u_0)\geq\eta$.
Then
\Eq{gf}{
  (g\oplus (n\odot f))(0)
  &=\sup_{\substack{u,v\in S\\ u+v=0}}\min\big(g(u),(n\odot f)(v)\big)\\
  &\geq \min\big(g(u_0),(n\odot f)(-u_0)\big)
  \geq \min\big(\eta,f(-\tfrac{u_0}n)\big).
}
There exists an open convex neighborhood $C$ of $0$ such that $f|_C=1$. Define the set $T_C$ by \eq{TA}, where the operation $*$ is given by \eq{+*}. Then, as we have seen it in the proof of \prp{3.2}, $T_C=G$ and $n*C\subseteq m*C$ holds for $n\leq m$. Therefore, there exists $n_0\in\N$ such that $-u_0\in n_0*C\subseteq n*C$ for all $n\geq n_0$. This means that for all $n\geq n_0$ there exists $c\in C$ and $w\in W$ such that $-u_0=n\cdot c+w$. Thus $-\frac{u_0}{n}=c+\frac{w}{n}$, which implies that $f(-\frac{u_0}{n})=f(c+\frac{w}{n})\geq f(c)=1$. Therefore, $f(-\frac{u_0}{n})=1$ for $n\geq n_0$. Combining this with the inequality \eq{gf}, we obtain
\Eq{*}{
  (g\oplus (n\odot f))(0)\geq \min\big(\eta,f(-\tfrac{u_0}n)\big)=\eta 
  \qquad(n\geq n_0).
}
Since $\eta$ was arbitrary, this inequality implies that $(g\oplus (n\odot f))(0)=1$, for all $n\geq n_0$. According to the first assertion of this theorem, this yields that $g\oplus (n\odot f)$ is a nonnegative element of the cornet $(F^1_W(S),\oplus,\odot,\leq)$. This proves that $f$ is Archimedean in $(F^1_W(S),\oplus,\odot,\leq)$.

(iii) Let $f$ be in $F_W^1(S)$ such that the support of $f$ is covered by the sum of a topologically bounded subset $D\subseteq S$ and $W$. Let $a\in\A$ be fixed. We need to show that there exists $n_0\in\N$ such that $f\leq n\odot a$ if $n\geq n_0$. 

By $a\in\A$, there exists an open convex neighborhood $C$ of $0$ such that $a|_C=1$. In view of the topological boundedness of $D$, we can find a number $n_0$ such that $D\subseteq n_0\cdot C$. By the convexity of $C$ and $0\in C$, this implies that $D\subseteq n\cdot C$. Consequently,
\Eq{*}{
  \supp(f)\subseteq D+W\subseteq n\cdot C+W.
}
Therefore, if $u\in\supp(f)$, then $\frac{u}{n}\in C+W$ for all $n\geq n_0$. The $W$-nondecreasingness of $a$ and $a|_C=1$ now yield that $a\big(\frac{u}{n}\big)=1$ for all $n\geq n_0$. This implies that $f(u)\leq a\big(\frac{u}{n}\big)=(n\odot a)(u)$ holds for all $u\in\supp(f)$ and $n\geq n_0$. This inequality is obvious for $u\in S\setminus\supp(f)$, i.e., if $f(u)=0$. Thus, we have proved that $f\leq n\odot a$ holds for all $n\geq n_0$, which shows that $f$ is $\A$-bounded.

(iv) in this part of the proof, we assume that $G$ is a locally convex topological group. Let $f\in F_W^1(S)$ be an upper semicontinuous element. To prove that $f$ is $\A$-closed, we need to show that if $g\in F_W^1(S)$ satisfies $g\leq f\oplus a$ for all $a\in\A$, then $g\leq f$. 

Let $x\in S$ be fixed and $\varepsilon>0$ be arbitrary. Then, by the upper semicontinuity of $f$ at $x$, there exists a neighborhood $U$ of $x$ such that $f(u)\leq f(x)+\varepsilon$ for all $u\in S\cap U$. Observing that $x-U$ is a neighborhood of $0$, the local convexity of $G$ implies that there exists a convex neighborhood $C$ of $0$ such that $C\subseteq x-U$, i.e., 
$x-C\subseteq U$. Let $a:=\chi_{C+W}$. Then $a$ is $W$-nondecreasing and $a|_C=1$, hence $a\in\A$. Therefore, we have $g\leq f\oplus a$, which implies
\Eq{*}{
  g(x)\leq (f\oplus a)(x)
  &=\sup_{\substack{u,v\in S\\ u+v=x}}\min\big(f(u),a(v)\big)
  =\sup_{\substack{u\in S,\, v\in C+W\\ u+v=x}} f(u)\\
  &=\sup_{\substack{u\in S\cap(x-C-W)}} f(u)
  \leq\sup_{\substack{u\in S\cap(U-W)}} f(u)
  \leq\sup_{\substack{u\in S\cap U}} f(u)\leq f(x)+\varepsilon.
}
Upon taking the limit $\varepsilon\to0$, the above inequality yields $f(x)\leq g(x)$, which was to be proved.

Finally, we prove the $\A$-continuity of the operation $\oplus$. Let $a\in \A$. Then there exists a $W$-invariant convex neighborhood $C$ of $0$ such that $a|_C=1$. Therefore, $\chi_C\leq a$. As we have seen it in the proof of \prp{3.2}, there exists a $W$-invariant convex neighborhood $B$ of $0$ such that $B+B\subseteq C$. Then 
\Eq{*}{
  \chi_B\oplus\chi_B=\Phi(B)\oplus\Phi(B)=\Phi(B+B)=\chi_{B+B}\leq\chi_C\leq a.
}
On the other hand $\chi_B\in\A$. This completes the proof of the $\A$-continuity of $\oplus$.
\end{proof}

\section{Main results}

The following result is the extension of the R{\aa}dström Cancellation Theorem to the setting of cornets.

\Thm{1}{
	Let $(X,+,*,\preceq)$ be a cornet and let $\A$ be a subsemigroup of $X_\prec$ such that the addition is $\A$-continuous. Let $x,y,z\in X$ such that $z$ is $\A$-bounded and $y$ is $\A$-closed and $m$-convex for some $m\geq2$. If
	\Eq{xyz}{
		x+z\preceq y+z
	}
	holds, then we have $x\preceq y$. 
}

\begin{proof}
	Let $x,y,z\in X$ such that \eq{xyz} holds. First, for all $n\in\N$, we show that
	\Eq{xyzn}{
		n\cdot x+z\preceq n\cdot y+z.
	}
	For $n=1$, the inequality is equivalent to \eq{xyz}. Assume that \eq{xyzn} holds for $n=k\in\N$. Then 
	\Eq{*}{
		(k+1)\cdot x+z
		=k\cdot x+x+z
		\preceq k\cdot y+x+z
		\preceq k\cdot y+y+z
		=(k+1)\cdot y+z,
	}
	which proves that the inequality \eq{xyzn} holds for $n=k+1$. This, by the principle of mathematical induction, completes the proof \eq{xyzn} for all $n\in\N$.

	The $m$-convexity of the element $y$ implies that $m\in C_y$, which is closed under multiplication by \lem{C}. Thus, for all $k\in\N$, $m^k\in C_y$. Using this and the inequality \eq{nx}, we obtain that
	\Eq{nxy}{
	  m^k*x\preceq m^k\cdot x \qquad\mbox{and}\qquad
	  m^k\cdot y= m^k*y.
	}
	Combining inequalities \eq{xyzn} and \eq{nxy}, it follows that for all $k\in\N$,
	\Eq{xyzk}{
		m^k*x+z\preceq m^k*y+z.
	}
	In the final step, assuming the $\A$-boundedness of $z$, we show that \eq{xyzk} implies $x\preceq y+a$ for all $a\in\A$. 
	
	Let $a\in\A$. Then, using that the addition is $\A$-continuous, we can find $b,c\in\A$ such that $b+c\preceq a$. Then there exists $n_1\in\N$ such that $0\preceq z+n*b$ holds for all $n_1\leq n$. The element $z$ is $\A$-bounded, thus we can find $n_2\in\N$ such that $z\preceq n*c$ holds for all $n_2\leq n$. By choosing $k$ so that $\max(n_1,n_2)\leq m^k$ is satisfied, it follows that
	\Eq{*}{
		0\preceq z+m^k*b\qquad\text{and}\qquad z\preceq m^k*c.
	}
	then we have
	\Eq{*}{
		m^k*x&\preceq m^k*x+z+m^k*b
		\preceq m^k*y+z+m^k*b \\
		&\preceq m^k*y+m^k*c+m^k*b=m^k*(y+c+b)
		\preceq m^k*(y+a).
	}
	This inequality implies that
	\Eq{*}{
	  x\preceq y+a
	}
	for all $a\in\A$. Now, using that $y$ is $\A$-closed, we can conclude that $x\preceq y$, which is what we wanted to prove.
\end{proof}

In what follows, we present several applications of the above theorem in the particular cornets described in \prp{1+}, \prp{2}, and \prp{3}.

\Cor{2}{
    Let $(G,+)$ be a locally convex topological abelian group such that there is no proper open subgroup of $G$. Let $W$ be a wedge and let $S$ be a subsemigroup of $G$ containing $W$. Let $P_W(S)$ denote the collection of\, $W$-invariant subsets of $S$ and define the operations $+$ and $*$ by \eq{+*}. Let $A,B,C\in P_W(S)$ such that $B$ is covered by the sum of a topologically bounded subset of $S$ and $W$ and $C$ is a topologically closed $m$-convex subset of $S$ for some $m\geq2$. Assume that
	\Eq{ABC}{
		A+B\subseteq C+B
	}
	holds. Then $A\subseteq C$. 
}

\begin{proof}
In view of \prp{3}, $(P_W(S),+,*,\subseteq)$ is a cornet and the $m$-convexity of $C$ implies that $C$ is an $m$-convex element of this cornet.

Let $\A$ denote the collection of those $W$-invariant subsets which contain an open convex neighborhood of $0$. Then, by assertion (ii) of \prp{3.2}, $\A$ is a subsemigroup of the Archimedean elements of $(P_W(S),+,*,\subseteq)$. By assertion (iii) of this proposition, we have that $B$ is $\A$-bounded and by the assertion (iv), we obtain that the element $C$ is $\A$-closed and the addition is $\A$-continuous. 

Therefore, we can apply \thm{1}, which shows that the inclusion \eq{ABC} implies $A\subseteq C$.
\end{proof}

\Cor{3}{
    Let $(G,+)$ be a locally convex topological abelian group such that there is no proper open subgroup of $G$. Let $W$ be a wedge and let $S$ be a uniquely divisible subsemigroup of $G$ containing $W$. Let the set $F^1_W(S)$ be defined by \eq{Fp} and define the operations $\oplus$ and $\odot$ by \eq{oper}. 	Let $f,g,h\in F^1_W(S)$ such that $\supp(h)$ is covered by the sum of a topologically bounded subset of $S$ and $W$ and $g$ is upper semicontinuous and $m$-quasiconcave for some $m\geq2$. Assume that
	\Eq{fgh}{
		f\oplus h\leq g\oplus h
	}
	holds. Then $f\leq g$. 
}

\begin{proof}
In view of \prp{3}, $(F^1_W(S),\oplus,\odot,\leq)$ is a cornet and the $m$-quasiconcavity of $g$ implies that $g$ is an $m$-convex element of this cornet.

Let $\A$ denote the collection of those $a\in F^1_W(S)$ for which there exists an open convex neighborhood $C$ of $0$ such that $a|_C=1$. Then, by assertion (ii) of \prp{3.3}, $\A$ is a subsemigroup of the Archimedean elements of $(F^1_W(S),\oplus,\odot,\leq)$. By assertion (iii) of this proposition, we have that $h$ is $\A$-bounded and by the assertion (iv), we obtain that the element $g$ is $\A$-closed and the addition $\oplus$ is $\A$-continuous. 

Therefore, we can apply \thm{1}, which shows that the inequality \eq{fgh} implies $f\leq g$.
\end{proof}

%\bibliography{publ,funcequ}

\begin{thebibliography}{10}

\bibitem{AghNou16}
M.~Aghajani and K.~Nourouzi, \emph{{On {H}ukuhara's differentiable iteration
  semigroups of linear set-valued functions}}, Aequationes Math. \textbf{90}
  (2016), no.~6, 1129–1145. \MR{3575583}

\bibitem{AghNouReg15}
M.~Aghajani, K.~Nourouzi, and D.~O'Regan, \emph{{The continuity of linear and
  sublinear correspondences defined on cones}}, Bull. Iranian Math. Soc.
  \textbf{41} (2015), no.~1, 43–55. \MR{3317175}

\bibitem{AzzBou15b}
D.~Azzam-Laouir and W.~Boukrouk, \emph{{A delay second-order set-valued
  differential equation with {H}ukuhara derivatives}}, Numer. Funct. Anal.
  Optim. \textbf{36} (2015), no.~6, 704–729. \MR{3349089}

\bibitem{AzzBou15a}
D.~Azzam-Laouir and W.~Boukrouk, \emph{{Second-order set-valued differential equations with boundary
  conditions}}, J. Fixed Point Theory Appl. \textbf{17} (2015), no.~1,
  99–121. \MR{3392984}

\bibitem{AzoGueMatMer10}
A.~Azócar, J.~A. Guerrero, J.~Matkowski, and N.~Merentes, \emph{{Uniformly
  continuous set-valued composition operators in the spaces of functions of
  bounded variation in the sense of {W}iener}}, Opuscula Math. \textbf{30}
  (2010), no.~1, 53–60. \MR{2591850}

\bibitem{BaiMosPop18}
A.~R. Baias, B.~Moşneguţu, and D.~Popa, \emph{{Set-valued solutions of a
  generalized quadratic functional equation}}, Results Math. \textbf{73}
  (2018), no.~4, Paper No. 129, 8. \MR{3851788}

\bibitem{BaiFar10}
R.~Baier and E.~Farkhi, \emph{{The directed subdifferential of {DC}
  functions}}, {Nonlinear analysis and optimization {II}. {O}ptimization},
  {Contemp. Math.}, vol. 514, Amer. Math. Soc., Providence, RI, 2010,
  p.~27–43. \MR{2668252}

\bibitem{BaiFarRos12a}
R.~Baier, E.~Farkhi, and V.~Roshchina, \emph{{The directed and {R}ubinov
  subdifferentials of quasidifferentiable functions, {P}art {I}: definition and
  examples}}, Nonlinear Anal. \textbf{75} (2012), no.~3, 1074–1088.
  \MR{2861321}

\bibitem{BaiFarRos12b}
R.~Baier, E.~Farkhi, and V.~Roshchina, \emph{{The directed and {R}ubinov subdifferentials of
  quasidifferentiable functions, {P}art {II}: calculus}}, Nonlinear Anal.
  \textbf{75} (2012), no.~3, 1058–1073. \MR{2861320}

\bibitem{BenSim16}
T.~Bendit and B.~Sims, \emph{{The structure of the normed lattice generated by
  the closed, bounded, convex subsets of a normed space}}, J. Nonlinear Convex
  Anal. \textbf{17} (2016), no.~6, 1069–1081. \MR{3540227}

\bibitem{BakGab09}
A.~Bąkowska and G.~Gabor, \emph{{Topological structure of solution sets to
  differential problems in {F}réchet spaces}}, Ann. Polon. Math. \textbf{95}
  (2009), no.~1, 17–36. \MR{2466011}

\bibitem{CheZho14}
Lixin Cheng and Yu~Zhou, \emph{{Approximation by {DC} functions and application
  to representation of a normed semigroup}}, J. Convex Anal. \textbf{21}
  (2014), no.~3, 651–661. \MR{3243811}

\bibitem{ChiTre13}
V.~V. Chistyakov and Y.~V. Tretyachenko, \emph{{A pointwise selection principle
  for maps of several variables via the total joint variation}}, J. Math. Anal.
  Appl. \textbf{402} (2013), no.~2, 648–659. \MR{3029178}

\bibitem{CorGajThi10}
R.~Correa, P.~Gajardo, and L.~Thibault, \emph{{Various {L}ipschitz-like
  properties for functions and sets. {I}. {D}irectional derivative and
  tangential characterizations}}, SIAM J. Optim. \textbf{20} (2010), no.~4,
  1766–1785. \MR{2600238}

\bibitem{CreKurRoc14}
G.~P. Crespi, D.~Kuroiwa, and M.~Rocca, \emph{{Convexity and global
  well-posedness in set-optimization}}, Taiwanese J. Math. \textbf{18} (2014),
  no.~6, 1897–1908. \MR{3284037}

\bibitem{DanMedMag11}
S.~Dancs, P.~Medvegyev, and Gy. Magyarkuti, \emph{{Normability via the
  convergence of closed and convex sets}}, J. Optim. Theory Appl. \textbf{150}
  (2011), no.~3, 675–682. \MR{2826296}

\bibitem{deBTom11}
F.~S. de~Blasi and L.~Tomassini, \emph{{On the strong law of large numbers in
  spaces of compact sets}}, J. Convex Anal. \textbf{18} (2011), no.~1,
  285–300. \MR{2777610}

\bibitem{Dol15}
M.~V. Dolgopolik, \emph{{Abstract convex approximations of nonsmooth
  functions}}, Optimization \textbf{64} (2015), no.~7, 1439–1469.
  \MR{3340642}

\bibitem{GayGeoMar16}
M.~Gaydu, M.~H. Geoffroy, and Y.~Marcelin, \emph{{Prederivatives of convex
  set-valued maps and applications to set optimization problems}}, J. Global
  Optim. \textbf{64} (2016), no.~1, 141–158. \MR{3437978}

\bibitem{GeoMar18}
M.~H. Geoffroy and Y.~Marcelin, \emph{{A concept of inner prederivative for
  set-valued mappings and its applications}}, ESAIM Control Optim. Calc. Var.
  \textbf{24} (2018), no.~3, 1059–1074. \MR{3877193}

\bibitem{Gon12}
Xiaobing Gong, \emph{{Convex solutions of the multi-valued iterative equation
  of order {$n$}}}, J. Inequal. Appl. (2012), 2012:258, 10. \MR{3017144}

\bibitem{GonSha08}
Z.~Gong and Y.~Shao, \emph{{Global existence and uniqueness of solutions for
  fuzzy differential equations under dissipative-type conditions}}, Comput.
  Math. Appl. \textbf{56} (2008), no.~10, 2716–2723. \MR{2460081}

\bibitem{GrzKucKucUrb14}
J.~Grzybowski, M.~Küçük, Y.~Küçük, and R.~Urbański,
  \emph{{Minkowski-{R}ådström-{H}örmander cone}}, Pac. J. Optim. \textbf{10}
  (2014), no.~4, 649–666. \MR{3275587}

\bibitem{GrzKucKucUrb15}
J.~Grzybowski, M.~Küçük, Y.~Küçük, and R.~Urbański, \emph{{On minimal representations by a family of sublinear
  functions}}, J. Global Optim. \textbf{61} (2015), no.~2, 279–289.
  \MR{3306087}

\bibitem{GrzPalPrzUrb12}
J.~Grzybowski, D.~Pallaschke, H.~Przybycień, and R.~Urbański,
  \emph{{Commutative semigroups with cancellation law: a representation
  theorem}}, Semigroup Forum \textbf{83} (2011), no.~3, 447–456. \MR{2860703}

\bibitem{GrzPalPrzUrb18}
J.~Grzybowski, D.~Pallaschke, H.~Przybycień, and R.~Urbański, \emph{{Reduced and minimally convex pairs of sets}}, J. Convex Anal.
  \textbf{25} (2018), no.~4, 1319–1334. \MR{3818568}

\bibitem{GrzPalUrb10a}
J.~Grzybowski, D.~Pallaschke, and R.~Urbański, \emph{{A note on the dual of
  the {M}inkowski-{R}ådström-{H}örmander lattice}}, Pac. J. Optim.
  \textbf{6} (2010), no.~2, 255–262. \MR{2668339}

\bibitem{GrzPalUrb16}
J.~Grzybowski, D.~Pallaschke, and R.~Urbański, \emph{{Reduced pairs of compact convex sets and ordered median
  functions}}, J. Optim. Theory Appl. \textbf{171} (2016), no.~2, 354–364.
  \MR{3557427}

\bibitem{GrzPalUrb19}
J.~Grzybowski, D.~Pallaschke, and R.~Urbański, \emph{{The formulas for the representation of functions of two
  variables as a difference of sublinear functions}}, Optimization \textbf{68}
  (2019), no.~10, 2055–2070. \MR{4003795}

\bibitem{GrzPrz15}
J.~Grzybowski and H.~Przybycień, \emph{{Completeness in
  {M}inkowski-{R}ådström-{H}örmander spaces}}, Optimization \textbf{64}
  (2015), no.~3, 495–503. \MR{3311974}

\bibitem{GrzPrz17}
J.~Grzybowski and H.~Przybycień, \emph{{Minimal representation in a quotient space over a lattice of
  unbounded closed convex sets}}, J. Convex Anal. \textbf{24} (2017), no.~2,
  695–705. \MR{3639284}

\bibitem{GrzPrzUrb13}
J.~Grzybowski, H.~Przybycień, and R.~Urbański, \emph{{Decomposition of
  {M}inkowski–{R}ådström–{H}örmander space to the direct sum of
  symmetric and asymmetric subspaces}}, Set-Valued Var. Anal. \textbf{21}
  (2013), no.~2, 201–216. \MR{3048244}

\bibitem{GrzUrb14}
J.~Grzybowski and R.~Urbański, \emph{{Dual space of the
  {M}inkowski-{R}ådström-{H}örmander space over {$\Bbb{R}^2$}}}, Funct.
  Approx. Comment. Math. \textbf{50} (2014), no.~1, [2013 on table of
  contents], 199–206. \MR{3189509}

\bibitem{HuaNin17}
Hui Huang and Jixian Ning, \emph{{Prederivatives of gamma paraconvex set-valued
  maps and {P}areto optimality conditions for set optimization problems}}, J.
  Inequal. Appl. (2017), Paper No. 243, 11. \MR{3708024}

\bibitem{Iva15}
G.~E. Ivanov, \emph{{Continuity and selections of the intersection operator
  applied to nonconvex sets}}, J. Convex Anal. \textbf{22} (2015), no.~4,
  939–962. \MR{3436695}

\bibitem{JouSil20}
A.~Jourani and F.~J. Silva, \emph{{Existence of {L}agrange multipliers under
  {G}âteaux differentiable data with applications to stochastic optimal
  control problems}}, SIAM J. Optim. \textbf{30} (2020), no.~1, 319–348.
  \MR{4060459}

\bibitem{Kho19}
H.~Khodaei, \emph{{Selections of generalized convex set-valued functions
  satisfying some inclusions}}, J. Math. Anal. Appl. \textbf{474} (2019),
  no.~2, 1104–1115. \MR{3926157}

\bibitem{Kur09}
D.~Kuroiwa, \emph{{On derivatives of set-valued maps and optimality conditions
  for set optimization}}, J. Nonlinear Convex Anal. \textbf{10} (2009), no.~1,
  41–50. \MR{2515285}

\bibitem{KurPopRoc15}
D.~Kuroiwa, N.~Popovici, and M.~Rocca, \emph{{A characterization of
  cone-convexity for set-valued functions by cone-quasiconvexity}}, Set-Valued
  Var. Anal. \textbf{23} (2015), no.~2, 295–304. \MR{3342735}

\bibitem{Kwi14}
G.~Kwiecińska, \emph{{On the {C}arathéodory superposition of multifunctions
  and an existence theorem}}, Math. Slovaca \textbf{64} (2014), no.~2,
  315–332. \MR{3201347}

\bibitem{LeiMerNikSan14}
H.~Leiva, N.~Merentes, K.~Nikodem, and J.~L. Sánchez, \emph{{Strongly convex
  set-valued maps}}, J. Global Optim. \textbf{57} (2013), no.~3, 695–705.
  \MR{3119375}

\bibitem{Mai12}
E.~Mainka-Niemczyk, \emph{{Some properties of set-valued sine families}},
  Opuscula Math. \textbf{32} (2012), no.~1, 159–170. \MR{2852477}

\bibitem{Mal15}
M.~T. Malinowski, \emph{{Set-valued and fuzzy stochastic differential equations
  in {M}-type 2 {B}anach spaces}}, Tohoku Math. J. (2) \textbf{67} (2015),
  no.~3, 349–381. \MR{3420550}

\bibitem{MirMah14}
A.~K. Mirmostafaee and M.~Mahdavi, \emph{{Approximately midconvex set-valued
  functions}}, Bull. Malays. Math. Sci. Soc. (2) \textbf{37} (2014), no.~2,
  525–530. \MR{3188055}

\bibitem{Orl17}
I.~V. Orlov, \emph{{On the embedding of a uniquely divisible {A}belian
  semigroup in a convex cone}}, Mat. Zametki \textbf{102} (2017), no.~3,
  396–404. \MR{3691704}

\bibitem{Pis11}
M.~Piszczek, \emph{{On multivalued iteration semigroups}}, Aequationes Math.
  \textbf{81} (2011), no.~1-2, 97–108. \MR{2773092 (2012a:39037)}

\bibitem{Pis13}
M.~Piszczek, \emph{{On selections of set-valued inclusions in a single variable
  with applications to several variables}}, Results Math. \textbf{64} (2013),
  no.~1-2, 1–12. \MR{3095123}

\bibitem{Pia09}
B.~Piątek, \emph{{On the continuity of the integrable multifunctions}},
  Opuscula Math. \textbf{29} (2009), no.~1, 81–88. \MR{2480434}

\bibitem{PloSkr14}
A.~V. Plotnikov and N.~V. Skripnik, \emph{{Conditions for the existence of
  local solutions of set-valued differential equations with generalized
  derivative}}, Ukrainian Math. J. \textbf{65} (2014), no.~10, 1498–1513,
  Translation of Ukraïn. Mat. Zh. {{\bf{6}}5} (2013), no. 10, 1350–1362.
  \MR{3215645}

\bibitem{Rad52b}
H.~Rådström, \emph{{An embedding theorem for spaces of convex sets}}, Proc.
  Amer. Math. Soc. \textbf{3} (1952), 165–169. \MR{0045938 (13,659c)}

\bibitem{Sik15}
J.~Sikorska, \emph{{Set-valued orthogonal additivity}}, Set-Valued Var. Anal.
  \textbf{23} (2015), no.~3, 547–557. \MR{3376895}

\bibitem{Sik16}
J.~Sikorska, \emph{{A singular behaviour of a set-valued approximate orthogonal
  additivity}}, Results Math. \textbf{70} (2016), no.~1-2, 163–172.
  \MR{3534999}

\bibitem{Sik19}
J.~Sikorska, \emph{{On a method of solving some functional equations for set-valued
  functions}}, Set-Valued Var. Anal. \textbf{27} (2019), no.~1, 295–304.
  \MR{3917368}

\bibitem{SmaSma12}
A.~Smajdor and W.~Smajdor, \emph{{Concave iteration semigroups of linear
  continuous set-valued functions}}, Cent. Eur. J. Math. \textbf{10} (2012),
  no.~6, 2272–2282. \MR{2983162}

\bibitem{Sma09}
W.~Smajdor, \emph{{On set-valued solutions of a functional equation of
  {D}rygas}}, Aequationes Math. \textbf{77} (2009), no.~1-2, 89–97.
  \MR{2495719 (2010a:39016)}

\bibitem{Sun17}
Yan Sun, \emph{{Asymptotic tests for interval-valued means}}, Statist. Probab.
  Lett. \textbf{121} (2017), 70–77. \MR{3575412}

\bibitem{Szc09}
J.~Szczawińska, \emph{{On some families of set-valued functions}}, Aequationes
  Math. \textbf{78} (2009), no.~1-2, 157–166. \MR{2552530 (2010i:54016)}

\bibitem{Szc13}
J.~Szczawińska, \emph{{On some equation for set-valued functions}}, Aequationes Math.
  \textbf{85} (2013), no.~3, 421–428. \MR{3063878}

\bibitem{VinNag12}
Cs. Vincze and A.~Nagy, \emph{{On the theory of generalized conics with
  applications in geometric tomography}}, J. Approx. Theory \textbf{164}
  (2012), no.~3, 371–390. \MR{2872524}

\bibitem{VinNag15}
Cs. Vincze and Á~Nagy, \emph{{Generalized conic functions of hv-convex planar
  sets: continuity properties and relations to {X}-rays}}, Aequationes Math.
  \textbf{89} (2015), no.~4, 1015–1030. \MR{3359690}

\bibitem{XuNikZha11}
B.~Xu, K.~Nikodem, and Zhang W., \emph{{On a multivalued iterative equation of
  order {$n$}}}, J. Convex Anal. \textbf{18} (2011), no.~3, 673–686.
  \MR{2858087 (2012i:39022)}

\end{thebibliography}
%\bibliographystyle{amsplain}
%\vfill\pagebreak

\providecommand{\bysame}{\leavevmode\hbox to3em{\hrulefill}\thinspace}
\providecommand{\MR}{\relax\ifhmode\unskip\space\fi MR }
\providecommand{\MRhref}[2]{%
  \href{http://www.ams.org/mathscinet-getitem?mr=#1}{#2}
}
\providecommand{\href}[2]{#2}

\end{document}